\documentclass[a4paper,10pt]{article}

\usepackage[T1]{fontenc}
\usepackage[english]{babel}
\usepackage{amssymb}
\usepackage{amsmath}
\usepackage{verbatim}
\usepackage{array}
\usepackage{lscape}
\usepackage[dvips]{color,graphicx}
\usepackage[arrow,curve,matrix]{xy}
\usepackage{mathrsfs}

\newtheorem{lemme}{Lemma}[section]
\newtheorem{prop}[lemme]{Proposition}
\newtheorem{propriété}[lemme]{Property}
\newtheorem{thm}[lemme]{Theorem}
\newtheorem{defi}[lemme]{Definition}
\newtheorem{cor}[lemme]{Corollary}

\newcommand{\dsum}{\displaystyle\sum}
\newcommand{\fl}{\longrightarrow}
\newcommand{\flto}{\longmapsto}
\newcommand{\C}{\mathbb{C}}

\newcommand{\R}{\mathbb{R}}
\newcommand{\Z}{\mathbb{Z}}
\newcommand{\N}{\mathbb{N}}
\newcommand{\eps}{\varepsilon}
\newcommand{\D}{\mathcal{D}}
\newcommand{\F}{\mathcal{F}}

\newcommand{\ol}{\stackrel{\mathbb{L}}{\otimes}}

\newcommand{\ts}{t^{s +1}}
\newcommand{\tsp}{\mathbf{t^{s+1}}}
\newcommand{\noyau}{\mathscr{T}}
\newcommand{\tnoyau}{\otimes_{r_{2}^{-1} \mathcal{O}_{\C^p / \Z^p}} \noyau }

\newenvironment{preuve}{\noindent\textsc{Proof}~:}{\hfill$\Box$\vspace*{0.8 cm}}

\title{\textbf{A commutation theorem on Mellin transform on sheaves and the functor of solutions \footnote{Preprint submitted to \emph{Annales de l'Institut Fourier} on October 5, 2006. This work was established at the laboratory Jean-Alexandre Dieudonné UMR CNRS 6621 (University of Nice Sophia-Antipolis, France).}
}}

\author{{\Large Hervé FABBRO}\\~\\
E-mail: Herve.Fabbro@ac-nice.fr}

\date{}

\begin{document}

\maketitle

\begin{center}
\begin{minipage}{0.75\textwidth}
\begin{center}
\textbf{Abstract}
\end{center}
We describe the complex of solutions of the algebraic
Mellin transform of a $\mathcal{D}$-module $\mathcal{M}$ in terms of
the solutions of $\mathcal{M}$. In order to do that, we define a
Mellin functor on sheaves. We show the Mellin transform of the
complex of rapid decay solutions of a regular holonomic
$\mathcal{D}$-module $\mathcal{M}$ is quasi-isomorphic to the
complex of solutions of the algebraic Mellin transform of
$\mathcal{M}$, the assumption of regularity not being necessary in the
one variable case.
\normalsize
\end{minipage}
\tableofcontents
\end{center}

\emph{\textbf{Acknowledgments} -- I offer heartfelt thanks to Michel Merle for his comments about the redaction of this preprint.}

					\section*{Introduction}

\indent Let us consider the affine space $\C^p$ $(p\geqslant 1)$ with coordinates $s_1,\dots,s_p$ and translations $\tau_j : s_j \flto s_j +1$ ($j$ from $1$ to $p$). We denote by $\mathbb{C}[s] \langle \tau, \tau^{-1} \rangle$ the non-commutative $\C$-algebra generated by the $s_j$, $\tau_j$, $\tau_j^{-1}$ and the relations $\tau_j s_j =(s_j+1) \tau_j$ ($j$ from $1$ to $p$). It is the finite difference operators algebra on $\mathbb{C}^p$, which is isomorphic to the Weil algebra $\mathrm{D}=\mathbb{C}[t,t^{-1}]\langle t\partial_t \rangle$ of polynomial differential operators 
on $(\C^*)^p$ through the correspondence $\tau_j \leftrightarrow t_j$ and $s_j \leftrightarrow -t_j\partial_{t_j}$.\\
The Mellin transform $\mathfrak{M}(\mathcal{M})$ of an algebraic $\mathcal{D}_{(\mathbb{C}^*)^p}$-module $\mathcal{M}$ is its global sections module $\mathrm{M}$ viewed as a $\mathbb{C}[s] \langle \tau, \tau^{-1} \rangle$, \emph{i.e.} equipped with the following action of $\mathbb{C}[s] \langle \tau, \tau^{-1} \rangle$: the action of $s_j$ is the one of $t_j\partial_{t_j}$ and that of $\tau_j$ is the one of the left multiplication by $t_j$ (F. Loeser and C. Sabbah studied this transformation in \cite{L-S}). Thus, we have a functor $\mathfrak{M} : \mathfrak{Mod}(\D_{(\mathbb{C}^*)^p})  \fl  \mathfrak{Mod}(\C[s] \langle \tau, \tau^{-1} \rangle)$.\\
Let $\C^p/\Z^p$ be the quotient of $\C^p$ by the group of the translations $\tau_j$, the $T_j=e^{-2i\pi s_j}$ the invariant functions and $\pi : \C^p \fl \C^p/\Z^p$ the canonical morphism. If $\mathcal{O}_{\C^p}$ is the sheaf of analytic functions on $\C^p$, we define the complex of solutions of a finite difference module $\mathbb{M}$ by $\mathcal{S}ol (\mathbb{M}) = \textbf{R}\mathcal{H}om_{\C[s]\langle \tau,\tau^{-1}\rangle} ( \mathbb{M} , \pi_* \mathcal{O}_{\C^p} )$, where $\C[s]\langle \tau,\tau^{-1}\rangle$ and $\mathbb{M}$ are viewed as constant sheaves. It is a complex of modules over the sheaf of analytic functions $\mathcal{O}_{\C^p/\Z^p}$ on $\C^p/\Z^p$. This is a definition similar to that of the complex of analytic solutions of a $\mathrm{D}$-module $\mathrm{M}$, \emph{i.e.} the complex $\mathcal{S}ol (\mathcal{M}) = \textbf{R}\mathcal{H}om_{\mathrm{D}} ( \mathrm{M} , \mathcal{O}_{(\C^*)^p} )$, where $\mathrm{D}$ and $\mathrm{M}$ are viewed as constant sheaves.\\
In 1992, C. Sabbah raised the following question:\\

\textbf{Question} (\cite{Sab}, Q$_2$ p.374) -- \emph{Can we define a Mellin transform functor on sheaves $\mathfrak{M}$ such that the following diagram is commutative?}
$$\xymatrix@R=1.5cm@C=2cm{\mathfrak{Mod}(\D_{(\C^*)^p})  \ar[r]^-{\mathcal{S}ol} \ar[d]^-{\mathfrak{M}} & D^b((\C^*)^p,\C) \ar@{-->}[d]^-{\mathfrak{M}} \\
\mathfrak{Mod}(\C[s]\langle \tau,\tau^{-1}\rangle) \ar[r]^-{\mathcal{S}ol} & D^b(\mathfrak{Mod}(\mathcal{O}_{\C^p/\Z^p})) }$$

Our way to deal with this problem is inspired by the work of B. Malgrange who solved a similar problem in the case of the Fourier transform (\cite{Mal-livre},\cite{Mal-Bourbaki}).\\

The classical Mellin transform being an integral transformation, it is natural to impose rapid decay conditions on the functions. The first step is to compactify the torus $(\C^*)^p$ into the annulus $(\overline{\C^*})^p$ which is identified with the real blowing-up of the divisor $(\mathbb{P}^1)^p \setminus (\C^*)^p$ in $(\mathbb{P}^1)^p$, and to use the sheaf $\mathcal{A}_{(\overline{\C^*})^p}^{<0}$ on $(\overline{\C^*})^p$ whose local sections are holomorphic on $(\C^*)^p$ and $\mathscr{C}^{\infty}$-flat on the boundary $(\overline{\C^*})^p \setminus (\C^*)^p$. Then, we define a rapid decay solutions functor by 
$$\mathcal{S}ol^{<0}~:~\mathcal{M} \flto \textbf{R}\mathcal{H}om_{\bar{\pi}^{-1} \D_{(\mathbb{P}^1)^p}^{an}} \big(\bar{\pi}^{-1} ( j_+ \mathcal{M})^{an}, \mathcal{A}_{(\overline{\C^*})^p}^{<0} \big)$$ 
where $j : (\C^*)^p \hookrightarrow  (\mathbb{P}^1)^p$ is the inclusion and $\bar{\pi}:(\overline{\C^*})^p \fl (\mathbb{P}^1)^p$ is the real blow-up described above.\\
This functor is more selective than the classical one. For example, the complexes of solutions of $t\dfrac{\partial}{\partial_t}$ and $t\dfrac{\partial}{\partial_t}+t$ are different for the functor $\mathcal{S}ol^{<0}$ but equal for the functor $\mathcal{S}ol$. This fact is very interesting because the constant function equal to $1$ and the exponential function obviously do not play the same role in this setup.\\  
The second step is the definition of the Mellin functor on sheaves. Following an idea of C. Sabbah (\cite{Sab}) dealing with Alexander modules of a complex of sheaves $\F$, we define the Mellin functor on sheaves $\mathfrak{M}$ by
$$\mathfrak{M}~:~ \F  \flto  \textbf{R}\Gamma \Big( (\overline{\C^*})^p , \F \ol_{\C}
\bar{\mathcal{L}} \Big)[p] \ol_{\C[T,T^{-1}]} \mathcal{O}_{\C^p / \Z^p}$$
where $\bar{\mathcal{L}}$ is the local system of rank $1$ free $\C[M_1,M_1^{-1}, \dots , M_p,M_p^{-1}]$-modules on $(\overline{\C^*})^p$ and the coordinates $T_1,\dots,T_p$ of $\C^p / \Z^p$ act on $\textbf{R}\Gamma \Big( (\overline{\C^*})^p , \F \ol_{\C} \bar{\mathcal{L}}\Big)$ by $M_1^{-1},\dots,M_p^{-1}$.\\

Then, our main \textbf{theorem \ref{th-mellin}} states that there exists a quasi-isomorphism
$$ \mathfrak{M}\Big(\mathcal{S}ol^{<0}(\mathcal{M})\Big) \cong \mathcal{S}ol \Big(\mathfrak{M}(\mathcal{M})\Big) $$
under suitable conditions for $\mathcal{M}$.

					\section{The main theorem}

\subsection{Rapid decay conditions on functions}\label{fast-dec}

Let us consider the compact annulus $\overline{\C^*}$ which is the real blow-up of $\mathbb{P}^1$ at $0$ and infinity.
We denote by $S_0$ and $S_\infty$ the circles $S^1$ over $0$ and infinity of $\overline{\C^*}$.\\
Then, we consider the following diagram in the category of real analytic manifolds with boundary (see a similar diagram in \cite{Da}):
$$\xymatrix{  &  (\overline{\C^*})^p  &  & \C^p / \Z^p \\
                     &  & (\overline{\C^*})^p  \times \C^p/\Z^p \ar[ul]_-{r_1} \ar[ur]^-{r_2} &  \\ 
                   & & (\overline{\C^*})^p  \times \C^p \ar[u]^-{\pi}  \ar[dd]^-{\bar{\pi}} &  \\
                  (\C^*)^p \times \C^p    \ar@{^{(}->}[urr]^-{k}  \ar@{^{(}->}[drr]^-{j} &  &  &  \\
                     &  & (\mathbb{P}^1)^p \times \C^p \ar[d]^-{\pi} & \\
                    &  & (\mathbb{P}^1)^p \times \C^p / \Z^p \ar[dl]^-{q_1} \ar[dr]_-{q_2}&\\
                   &   (\mathbb{P}^1)^p &  &\C^p / \Z^p}$$
Let us recall that $\pi : \C^p \fl \C^p / \Z^p$ is the quotient by actions of $\tau_1, \dots ,\tau_p$. We write $\pi$, $\bar{\pi}$, $k$ and $j$ instead of $id \times \pi$, $\bar{\pi} \times id_{\C^p}$, $k \times id_{\C^p}$ and $j \times id_{\C^p}$.\\

Let us recall some notations (for example, see \cite{Sab-Stokes-2} chapter II.1). We denote by $Y$ any analytic variety.\\ 
Let $\mathcal{C}^{\infty}_{(\mathbb{P}^1)^p \times Y}$ be the sheaf of functions infinitly differentiable with respect to $(\mathbb{P}^1)^p$ and holomorphic with respect to $Y$, and $\mathcal{C}^{<0}_{(\mathbb{P}^1)^p \times Y}$ be the subsheaf of $\mathcal{C}^{\infty}_{(\mathbb{P}^1)^p \times Y}$ of flat functions on $((\mathbb{P}^1)^p \setminus (\C^*)^p) \times Y$. Likewise, let $\mathcal{C}^{\infty}_{(\overline{\C^*})^p \times Y}$ be the sheaf of functions infinitly differentiable with respect to $(\overline{\C^*})^p$ and holomorphic with respect to $Y$, and $\mathcal{C}^{<0}_{(\overline{\C^*})^p \times Y}$ be the subsheaf of $\mathcal{C}^{\infty}_{(\overline{\C^*})^p \times Y}$ of flat functions on $((\overline{\C^*})^p \setminus (\C^*)^p) \times Y$.\\
We denote by
$$\mathcal{P}^{0,\bullet}_{(\mathbb{P}^1)^p \times Y} = \Big(\mathcal{C}^{<0}_{(\mathbb{P}^1)^p \times Y} \otimes_{\mathcal{C}^{\infty}_{(\mathbb{P}^1)^p \times Y}} \mathcal{C}^{0,\bullet}_{(\mathbb{P}^1)^p \times Y},\bar{\partial}\Big)$$
and 
$$\mathcal{P}^{0,\bullet}_{(\overline{\C^*})^p \times Y} = \Big(\mathcal{C}^{<0}_{(\overline{\C^*})^p \times Y} \otimes_{\mathcal{C}^{\infty}_{(\overline{\C^*})^p \times Y}} \mathcal{C}^{0,\bullet}_{(\overline{\C^*})^p \times Y},\bar{\partial}\Big)$$
the flat Dolbeault complexes, where $\mathcal{C}^{k,l}_{(\mathbb{P}^1)^p \times Y}$ and $\mathcal{C}^{k,l}_{(\overline{\C^*})^p \times Y}$ 
are the sheaves of infinitly differentiable $(k,l)$-forms which are holomorphic with respect to $Y$.\\

We have $\textbf{R}\bar{\pi}_*\mathcal{C}^{<0}_{(\overline{\C^*})^p  \times Y} = \bar{\pi}_*\mathcal{C}^{<0}_{(\overline{\C^*})^p  \times Y} = \mathcal{C}^{<0}_{(\mathbb{P}^1)^p \times Y}$ 
(see \cite{Mal-ideals}).\\ 
Thus, for all $(k,l)$, we have $\textbf{R}\bar{\pi}_*\mathcal{P}^{k,l}_{(\overline{\C^*})^p \times Y}=\bar{\pi}_*\mathcal{P}^{k,l}_{(\overline{\C^*})^p \times Y}=\mathcal{P}^{k,l}_{(\mathbb{P}^1)^p \times Y}$.\\

Then, we define
$$\xymatrix@C=1.5cm{\mathcal{A}_{(\overline{\C^*})^p  \times Y}^{<0} = \mathrm{Ker} \Big(  \mathcal{P}^{0,0}_{(\overline{\C^*})^p \times Y} \ar[r]^-{\bar{\partial}} & \mathcal{P}^{0,1}_{(\overline{\C^*})^p \times Y}}\Big)$$
It is the subsheaf of $\mathcal{C}^{<0}_{(\overline{\C^*})^p  \times Y}$ whose local sections are holomorphic on $(\C^*)^p\times Y$, namely the sheaf on $(\overline{\C^*})^p  \times Y$ made up of holomorphic functions with rapid decay condition regarding each $t_j$ in angular sectors along the circles $S^1$ at $0$ and at infinity, uniformly relative to the $s_j$ on any compact.\\
Thus, we get the resolution $\Big( \mathcal{P}^{0,\bullet}_{(\overline{\C^*})^p \times Y}, \bar{\partial}\Big)$ of $\mathcal{A}_{(\overline{\C^*})^p  \times Y}^{<0}$ in the category of $\mathrm{D}$-modules (\cite{Maj}). The morphism $\pi$ being a covering, we also have a quasi-isomorphim of complexes of $\mathrm{D}$-modules
$$\pi_*\mathcal{A}_{(\overline{\C^*})^p  \times \C^p}^{<0} \cong \Big( \pi_*\mathcal{P}^{0,\bullet}_{(\overline{\C^*})^p \times \C^p}, \bar{\partial}\Big)$$

In the continuation, the space $Y$ will be $\C^p$, $\C^p/\Z^p$ or a point, and in this case, we will identify $(\overline{\C^*})^p  \times \{*\}$ with $(\overline{\C^*})^p$.

\subsection{The functor of rapid decay solutions}

\begin{defi}
The functor of rapid decay holomorphic solutions of an algebraic $\D_{(\C^*)^p}$-module is 
$$\begin{array}{rrcl}
\mathcal{S}ol^{<0} : & \mathfrak{Mod}(\D_{(\C^*)^p}) & \fl & D^b((\overline{\C^*})^p,\C) \\
& \mathcal{M} & \flto & \textbf{R}\mathcal{H}om_{\bar{\pi}^{-1} \D_{(\mathbb{P}^1)^p}^{an}} \Big(\bar{\pi}^{-1} ( j_+ \mathcal{M})^{an}, \mathcal{A}_{(\overline{\C^*})^p}^{<0} \Big)
\end{array}$$
where $( j_+ \mathcal{M})^{an}$ is the analytised of $j_+ \mathcal{M}$.
\end{defi}

The morphism $j$ being the inclusion of an open affine space, we get $j_+ \mathcal{M} = \textbf{R}j_*\mathcal{M}=j_*\mathcal{M}$ (for example, see \cite{Bo}). Thus, we get
$$\mathcal{S}ol^{<0}(\mathcal{M}) = \textbf{R}\mathcal{H}om_{\bar{\pi}^{-1} \D_{(\mathbb{P}^1)^p}^{an}} \Big(\bar{\pi}^{-1} ( j_* \mathcal{M})^{an}, \mathcal{A}_{(\overline{\C^*})^p}^{<0} \Big)$$

\begin{lemme}\label{sol-dec-rapide}
For any algebraic $\D_{(\C^*)^p}$-module $\mathcal{M}$, whose module of global sections is denoted by $\mathrm{M}$, we have
$$\mathcal{S}ol^{<0}(\mathcal{M})=\textbf{R}\mathcal{H}om_{\mathrm{D}} ( \mathrm{M} , \mathcal{A}_{(\overline{\C^*})^p}^{<0} )$$
where $\mathrm{D}$ and $\mathrm{M}$ are viewed as constant sheaves.
\end{lemme}

\begin{preuve}
The module $\mathcal{M}$ being coherent on the affine space $(\C^*)^p$, we have easily that 
$$\bar{\pi}^{-1} (j_*\mathcal{M})^{an} = \bar{\pi}^{-1} (j_*\D_{(\C^*)^p})^{an}
\otimes_{\mathrm{D}_{(\overline{\C^*})^p}} \mathrm{M}_{(\overline{\C^*})^p}$$
Thanks to the adjonction of the tensor product and the functor $\textbf{R}\mathcal{H}om$, we have
\begin{eqnarray*}
\mathcal{S}ol^{<0}(\mathcal{M}) & =  & \textbf{R}\mathcal{H}om_{\mathrm{D}} \Big( \mathrm{M} , \textbf{R}\mathcal{H}om_{\bar{\pi}^{-1} \D_{(\mathbb{P}^1)^p}^{an}} (\bar{\pi}^{-1} (j_*\D_{(\C^*)^p})^{an},\mathcal{A}_{(\overline{\C^*})^p}^{<0} )\Big)
\end{eqnarray*}
The sheaf $(j_*\D_{(\C^*)^p})^{an}$ is a locally free $\D_{(\mathbb{P}^1)^p}^{an}$-module of rank $1$: in a neighbourhood of a point $t^0=(0,\dots,0,\infty,\dots,\infty,t_{r+1},\dots,t_p) \in \{0\}^l \times\{\infty\}^{r-l}\times (\C^*)^{p-r}\subset (\mathbb{P}^1)^p$, we have  
$$(j_*\D_{(\C^*)^p})^{an}_{t^0} = \D_{(\mathbb{P}^1)^p,t^0}^{an} \Big[\dfrac{1}{t_1},\dots,\dfrac{1}{t_l},t_{l+1},\dots,t_r\Big] = \D_{(\mathbb{P}^1)^p,t^0}^{an} \dfrac{t_{l+1} \cdots t_r}{t_1 \cdots t_l}$$
and there exists a morphism 
$$\D_{(\mathbb{P}^1)^p,t^0}^{an} \stackrel{u}{\fl} \D_{(\mathbb{P}^1)^p,t^0}^{an} \dfrac{t_{l+1} \cdots t_r}{t_1 \cdots t_l}$$
defined by $u(1)=1=\dfrac{t_1 \cdots t_l}{t_{l+1} \cdots t_r}.\dfrac{t_{l+1} \cdots t_r}{t_1 \cdots t_l}$.\\
In a neighbourhood of $\theta^0=(\theta_1,\dots,\theta_l,\theta_{l+1},\dots,\theta_{r},t_{r+1},\dots,t_p) \in \bar{\pi}^{-1}(t^0)$, proving that
$$\textbf{R}\mathcal{H}om_{\bar{\pi}^{-1} \D_{(\mathbb{P}^1)^p}^{an}} (\bar{\pi}^{-1} (j_*\D_{(\C^*)^p})^{an},\mathcal{A}_{(\overline{\C^*})^p}^{<0} ) \stackrel{u^*}{\fl} \mathcal{A}_{(\overline{\C^*})^p}^{<0}$$
is an isomorphism amounts to showing that
$$\mathcal{H}om_{\bar{\pi}^{-1} \D_{(\mathbb{P}^1)^p}^{an}} (\bar{\pi}^{-1} (j_*\D_{(\C^*)^p})^{an},\mathcal{A}_{(\overline{\C^*})^p}^{<0} ) \stackrel{u^*}{\fl} \mathcal{A}_{(\overline{\C^*})^p}^{<0}$$
is an isomorphism.\\
Now, $u^*$ being defined by $u^* : \varphi \flto \varphi_\circ u (1) =\varphi(1)$ and the morphism $\varphi$ being defined by $\varphi(1)$, $u^*$ is an isomorphism.\\
Finally, let us state precisely the structure of $\mathrm{D}$-module :
$$u^*(t_j \partial_{t_j}.\varphi)=(t_j \partial_{t_j}.\varphi)(1)=\varphi(1.t_j \partial_{t_j})=\varphi(t_j \partial_{t_j})=t_j \partial_{t_j}\varphi(1)=t_j \partial_{t_j}u^*(\varphi)$$
$$u^*(t_j .\varphi)=(t_j .\varphi)(1)=\varphi(1.t_j)=\varphi(t_j)=t_j \varphi(1)=t_j u^*(\varphi)$$
\end{preuve}

\subsection{Definition of the Mellin transform on sheaves}\label{def-mellin-sheaves}

Such a transformation was already studied by O. Gabber and F. Loeser in the $\ell$-adic case (\cite{G-L},\cite{Loe}).
C. Sabbah also announces a definition in the complex case (\cite{Sab}). 
In definition \ref{def-mellin-faisceautique}, we give an adaptation to our case.
 
Let us denote by $\bar{q}:\widetilde{(\overline{\C^*})^p} \fl (\overline{\C^*})^p$ the universal covering of $(\overline{\C^*})^p$, and
$$\C[M,M^{-1}]=\C[M_1,M_1^{-1}, \dots , M_p,M_p^{-1}]$$
$$\C[T,T^{-1}]=\C[T_1,T_1^{-1}, \dots , T_p,T_p^{-1}]$$
Let $\bar{\mathcal{L}}=\bar{q}_!\bar{q}^{-1}\underline{\C}_{(\overline{\C^*})^p}$ be the natural local system of rank $1$ free $\C[M,M^{-1}]$-modules on $(\overline{\C^*})^p$, where multiplications by the $M_j$ are  monodromies around the coordinates hyperplans.\\
 
\begin{defi}\label{def-mellin-faisceautique}
The Mellin transform functor on sheaves is
$$\begin{array}{rrcl}
\mathfrak{M} : & D^b((\overline{\C^*})^p,\C) & \fl & D^b(\mathfrak{Mod}(\mathcal{O}_{\C^p / \Z^p})) \\
& \F & \flto & \textbf{R}\Gamma \Big( (\overline{\C^*})^p , \F \ol_{\C}
\bar{\mathcal{L}} \Big)[p] \ol_{\C[T,T^{-1}]} \mathcal{O}_{\C^p / \Z^p}
\end{array}$$
where $T_j$ acts on $\textbf{R}\Gamma \Big( (\overline{\C^*})^p , \F \ol_{\C} \bar{\mathcal{L}}\Big)$ by
$M_j^{-1}$. 
\end{defi}

If $\F$ is an object of $D^b((\overline{\C^*})^p,\C)$ with $\R$-constructible cohomology sheaves, then $\mathfrak{M}(\F)$ is a complex of $\mathcal{O}_{\C^p / \Z^p}$-modules with coherent cohomology.

\subsection{Statement of the main theorem}

The commutation theorem we want to prove is the following:

\begin{thm}\label{th-mellin}
For any coherent algebraic $\mathcal{D}_{(\C^*)^p}$-module $\mathcal{M}$, such that 
$\mathcal{S}ol^{<0}(\mathcal{M})\big\vert_{(\C^*)^p}$ has $\C$-constructible cohomology and $\mathcal{S}ol^{<0}(\mathcal{M})\big\vert_{(\overline{\C^*})^p \setminus (\C^*)^p}$ has $\R$-constructible cohomology,
there exists a canonical isomorphism in $D^b(\mathfrak{Mod}(\mathcal{O}_{\C^p / \Z^p}))$
$$ \mathfrak{M}\Big(\mathcal{S}ol^{<0}(\mathcal{M})\Big) \cong \mathcal{S}ol \Big(\mathfrak{M}(\mathcal{M})\Big) $$
\end{thm}

A classical result in the theory of $\D$-modules sets that, if $\mathcal{M}$ is an holonomic $\D_{(\C^*)^p}$-module which is regular at $(\mathbb{P}^1)^p \setminus (\C^*)^p$, then the complex $\mathcal{S}ol^{<0}(\mathcal{M})$ has $\C$-constructible cohomology sheaves.\\
In the one variable case, if we abstain from the regularity, we know the following result (\cite{Mal-livre},\cite{Was}): for any holonomic $\D_{\C^*}$-module, the restriction of the complex $\mathcal{S}ol^{<0}(\mathcal{M})$ to $\C^*$ has $\C$-constructible cohomology, and its restriction to $\overline{\C^*} \setminus \C^*$ has $\R$-constructible cohomology.\\
Thanks to these results, we get the following corollary:

\begin{cor}\label{cor-au-th-mellin}
If one of the two following assumptions holds 
\begin{itemize}
\item $\mathcal{M}$ is an holonomic $\mathcal{D}_{(\C^*)^p}$-module which is regular at $(\mathbb{P}^1)^p \setminus (\C^*)^p$
\item for $p=1$, $\mathcal{M}$ is an holonomic $\mathcal{D}_{\C^*}$-module
\end{itemize}
then there exists a canonical isomorphism in $D^b(\mathfrak{Mod}(\mathcal{O}_{\C^p / \Z^p}))$
$$ \mathfrak{M}\Big(\mathcal{S}ol^{<0}(\mathcal{M})\Big) \cong \mathcal{S}ol \Big(\mathfrak{M}(\mathcal{M})\Big) $$
\end{cor}

There also exists a conjecture stating that, if $\mathcal{M}$ is an holonomic $\D_{(\C^*)^p}$-module, then the complex $\mathcal{S}ol^{<0}(\mathcal{M})$ has $\C$-constructible cohomology sheaves for any integer $p\geqslant 1$. It is well known on $(\C^*)^p$ (\cite{Kas},\cite{M-N}) but not on the divisor at infinity. So, we can formulate the conjecture:\\

\textbf{Conjecture} -- \emph{For any holonomic coherent algebraic $\mathcal{D}_{(\C^*)^p}$-module $\mathcal{M}$, there exists a canonical isomorphism in $D^b(\mathfrak{Mod}(\mathcal{O}_{\C^p / \Z^p}))$}
$$ \mathfrak{M}\Big(\mathcal{S}ol^{<0}(\mathcal{M})\Big) \cong \mathcal{S}ol \Big(\mathfrak{M}(\mathcal{M})\Big) $$

					\section{A finiteness theorem on solutions}

The aim of this section is to prove our second theorem which is useful in the proof of theorem \ref{th-mellin} but also interesting by itself .

\begin{thm}\label{finitude-solutions}
Consider a coherent algebraic $\D_{(\C^*)^p}$-module $\mathcal{M}$. We denote its module of global sections by $\mathrm{M}$.
Then, if 
$\textbf{R}\mathcal{H}om_{\mathrm{D}} ( \mathrm{M} , 
\mathcal{A}_{(\overline{\C^*})^p}^{<0})$ has $\C$-constructible cohomology on $(\C^*)^p$ and $\R$-constructible cohomology on $(\overline{\C^*})^p$,
we have a canonical isomorphism in $D^b(\mathfrak{Mod}(\mathcal{O}_{\C^p / \Z^p}))$
$$\textbf{R}\mathcal{H}om_{\mathrm{D}} ( \mathrm{M} ,
{r_1}^{-1}\mathcal{A}_{(\overline{\C^*})^p}^{<0}) \ol_{\C} r_{2}^{-1}\mathcal{O}_{\C^p / \Z^p} 
\cong \textbf{R}\mathcal{H}om_{\mathrm{D}} ( \mathrm{M} ,
\mathcal{A}_{(\overline{\C^*})^p \times \C^p/\Z^p}^{<0})$$
where $\mathrm{D}$ and $\mathrm{M}$ are viewed as constant sheaves.
\end{thm}

Let us prove this result.
There exists a canonical morphism\\
$\textbf{R}\mathcal{H}om_{\mathrm{D}} ( \mathrm{M} ,{r_1}^{-1}\mathcal{A}_{(\overline{\C^*})^p}^{<0}) \ol_{\C} r_{2}^{-1}\mathcal{O}_{\C^p / \Z^p}$

\hfill $\fl \textbf{R}\mathcal{H}om_{\mathrm{D}} ( \mathrm{M} ,{r_1}^{-1}\mathcal{A}_{(\overline{\C^*})^p}^{<0} \otimes_{\C} r_{2}^{-1}\mathcal{O}_{\C^p / \Z^p})$\\
which is an isomorphism because $\mathrm{M}$ is a finitely presented $\mathrm{D}$-module.\\

It is now sufficient to prove that the canonical morphism
$$\textbf{R}\mathcal{H}om_{\mathrm{D}}(\mathrm{M},{r_1}^{-1}\mathcal{A}_{(\overline{\C^*})^p}^{<0}\otimes_{\C} r_{2}^{-1}\mathcal{O}_{\C^p / \Z^p}) \stackrel{\Phi}{\fl} \textbf{R}\mathcal{H}om_{\mathrm{D}} ( \mathrm{M} , \mathcal{A}_{(\overline{\C^*})^p \times \C^p/\Z^p}^{<0})$$
is a quasi-isomorphism when the hypothesis of the theorem \ref{finitude-solutions} is satisfied.\\
On a free resolution of the $\mathrm{D}$-module $\mathrm{M}$, the stalk at a point $x^0=(\theta_I^0,\theta_J^0,t^0,T^0) \in (S_0)^I \times (S_\infty)^J \times (\C^*)^{\{1,\dots,p\}\setminus I\cup J}\times \C^p/\Z^p$ of the previous morphism is written
$$\xymatrix{0 \ar[r] & (\mathcal{A}_{(\overline{\C^*})^p}^{<0}\boxtimes_{\C} \mathcal{O}_{\C^p / \Z^p})_{x^0}^{m_0} \ar[d]^-{\Phi_0} \ar[r]^-{\Psi_0} & \dots \ar[r]^-{\Psi_{l-1}} & (\mathcal{A}_{(\overline{\C^*})^p}^{<0}\boxtimes_{\C} \mathcal{O}_{\C^p / \Z^p})_{x^0}^{m_l} \ar[r]^-{\Psi_{l}} \ar[d]^-{\Phi_l}& \dots  \\
0 \ar[r] & (\mathcal{A}_{(\overline{\C^*})^p \times \C^p/\Z^p,x^0}^{<0})^{m_0} \ar[r]^-{\Psi_0} & \dots \ar[r]^-{\Psi_{l-1}} & (\mathcal{A}_{(\overline{\C^*})^p \times \C^p/\Z^p,x^0}^{<0})^{m_l} \ar[r]^-{\Psi_{l}} & \dots}$$
where the $\Psi_j$ are matrices with entries in $\mathrm{D}$.\\

The complex $\textbf{R}\mathcal{H}om_{\mathrm{D}}(\mathrm{M},
\mathcal{A}_{(\overline{\C^*})^p}^{<0})$ has constructible cohomology sheaves. Then, the stalks of the latter are finite dimensional complex vector spaces.\\

\textbf{Let us prove that $\Phi_0$ induces an isomorphism on cohomology}\\

Let us recall the following notations:
$$T^\alpha=T_1^{\alpha_1} \cdots T_p^{\alpha_p}~\textrm{and}~(T-T^0)^\alpha=(T_1-T^0_1)^{\alpha_1} \cdots (T_p-T^0_p)^{\alpha_p}$$
$$t^N = t_1^{N_1} \cdots t_p^{N_p}$$
for all $\alpha=(\alpha_1,\dots,\alpha_p) \in \N^p$ and for all $N=(N_1,\dots,N_p) \in \N^p$.\\
For all $(\theta,\eps)=(\theta_1,\dots,\theta_q,\eps_1,\dots,\eps_q) \in (\R^q)^2$, let us denote by $D_{0}(\theta,\eps)$ the polysector $\displaystyle\prod_{i=1}^q \{ |t_i|<\eps_i~,~|\mathrm{Arg}(t_i)-\theta_i|<\eps_i \}$ and $D_{\infty}(\theta,\eps)$ the polysector $\displaystyle\prod_{i=1}^q \{ |t_i|>\dfrac{1}{\eps_i}~,~|\mathrm{Arg}(t_i)-\theta_i|<\eps_i \}$.\\

The morphism $\Phi_0$ is the inclusion. Let us show that $g = (g_1,\dots,g_{m_0}) \in (\mathcal{A}_{(\overline{\C^*})^p \times \C^p/\Z^p,x^0}^{<0})^{m_0}$ such that $\Psi_0 (g)=0$ is an element of 
$(\mathcal{A}_{(\overline{\C^*})^p}^{<0}\boxtimes_{\C} \mathcal{O}_{\C^p / \Z^p})_{x^0}^{m_0}$.\\

All sections $f \in \mathcal{A}_{(\overline{\C^*})^p \times \C^p/\Z^p,x^0}^{<0}$, where $x^0=(\theta_I^0,\theta_J^0,t^0,T^0) \in (S_0)^I \times (S_\infty)^J \times (\C^*)^{\{1,\dots,p\}\setminus I\cup J}\times \C^p/\Z^p$, are holomorphic in an open set $D_{0}(\theta_I^0,\eps_i)  \times D_{\infty}(\theta_J^0,\eps_j) \times W^0 \times V^0$, $W^0$ being an open neighbourhood of $t^0$ and $V^0$ an open neighbourhood of a compact polydisc  $\overline{D(T^0,R)}=\overline{D(T^0_1,R_1)} \times \cdots \times \overline{D(T^0_p,R_p)}$ with center $T^0$ and polyradius $R \in (\R^*_+)^p$.

\begin{lemme}\label{developpement}
All function $f \in \mathcal{A}_{(\overline{\C^*})^p \times \C^p/\Z^p,x^0}^{<0}$ defined in an open set described above expands into 
$$f = \sum_{\alpha \in \N^p} u_{\alpha}.(T - T^0)^\alpha~~\textrm{where}~~u_\alpha \in \mathcal{A}_{(\overline{\C^*})^p,(\theta_I^0,\theta_J^0,t^0)}^{<0}$$ 
and $\vert T_j - T^0_j \vert < R_j$.\\
This series converge normally, namely
$$\forall K\subset W^0~\textrm{compact},~\sum_{\alpha \in \N^p} \Vert u_{\alpha} \Vert_\infty^{K,\eps}~\rho^\alpha < +\infty$$
for $\rho=(\rho_1,\dots,\rho_p) \in (\R^*_+)^p$ such that $\rho_j < R_j$ for all $j=1,\dots,p$, and where $\Vert - \Vert_\infty^{K,\eps}$ is the norm of the uniform convergence in $D_{0}(\theta_I^0,\eps_i)  \times D_{\infty}(\theta_J^0,\eps_j) \times K$.
\end{lemme}

\begin{preuve}
thanks to Cauchy formula, we have for all $(t,T) \in D_{0}(\theta_I^0,\eps_i)  \times D_{\infty}(\theta_J^0,\eps_j) \times W^0 \times D(T^0,R)$:
$$f(t,T) = \Big(\frac{1}{2i\pi}\Big)^p \int \cdots \int_{ \{ \vert \xi_j - T^0_j \vert = R_j \} } \frac{f(t,\xi)}{(\xi_1 - T_1) \cdots (\xi_p - T_p)} d\xi_1 \cdots d\xi_p$$ 
In writing
$$\frac{1}{(\xi_1 - T_1) \cdots (\xi_p - T_p)} = \sum_{k_1,\dots,k_p \geqslant 0} \frac{(T_1-T^0_1)^{k_1} \cdots (T_p-T^0_p)^{k_p}}{(\xi_1-T^0_1)^{k_1 +1} \cdots (\xi_p-T^0_p)^{k_p +1}}$$
and noting
$$u_{k_1,\dots,k_p} (t) = \Big(\frac{1}{2i\pi}\Big)^p  \int \cdots \int_{ \{ \vert \xi_j - T^0_j \vert = R_j \} }  \frac{f(t,\xi)}{(\xi_1-T^0_1)^{k_1 +1} \cdots (\xi_p-T^0_p)^{k_p +1}} d\xi_1 \cdots d\xi_p $$
we get
$$f(t,T) = \sum_{k_1,\dots,k_p \geqslant 0} u_{k_1,\dots,k_p} (t) (T_1-T^0_1)^{k_1} \cdots (T_p-T^0_p)^{k_p}$$
By domination, the $u_{k_1,\dots,k_p}$ are holomorphic.\\
To get the rapid decay condition of the $u_{k_1,\dots,k_p}$, we use that of $f$. Let $K$ be a compact contained in $W^0$. For all $t\in D_{0}(\theta_I^0,\eps_i)  \times D_{\infty}(\theta_J^0,\eps_j) \times K$:  
\begin{eqnarray*}
\vert u_{k_1,\dots,k_p}(t) \vert & \leqslant & \frac{1}{(2\pi)^p} \Bigg\vert \int \cdots \int_{ \{ \vert \xi_j - T^0_j \vert = R_j \} }  \frac{\vert f(t,\xi) \vert}{\vert\xi_1-T^0_1\vert^{k_1 +1} \cdots \vert\xi_p-T^0_p\vert^{k_p +1}} d\xi_1 \cdots d\xi_p \Bigg\vert \\
 & \leqslant & \frac{1}{(2\pi)^p} \Bigg\vert \int \cdots \int_{ \{ \vert \xi_j - T^0_j \vert = R_j \} }  \frac{\vert f(t,\xi) \vert}{R_1^{k_1 +1} \cdots R_p^{k_p +1}} d\xi_1 \cdots d\xi_p \Bigg\vert \\
 & \leqslant & \frac{1}{(2\pi)^p} \Bigg\vert \int \cdots \int_{ \{ \vert \xi_j - T^0_j \vert = R_j \} }  \frac{C_{N} \frac{\prod_{i \in I} |t_i|^{N_i}}{ \prod_{j \in J} |t_j|^{N_j}} }{R_1^{k_1 +1} \cdots R_p^{k_p +1}} d\xi_1 \cdots d\xi_p \Bigg\vert \\
 & \leqslant &  \frac{C_{N}}{R_1^{k_1} \cdots R_p^{k_p}}  \frac{\prod_{i \in I} |t_i|^{N_i}}{ \prod_{j \in J} |t_j|^{N_j}}
\end{eqnarray*}
for all $N=(N_k)_{k \in I \cup J} \in \N^{I \cup J}$, the constant $C_N >0$ depending on $N$.\\
The convergence of the series results from this. For all $\rho\in (\R^*_+)^p$ such that $\rho_j < R_j$ for all $j=1,\dots,p$, we get:
\begin{eqnarray*}
\forall K\subset W^0~\textrm{compact},~\sum_{\alpha \in \N^p} \Vert u_{\alpha} \Vert_\infty^{K,\eps}~\rho^\alpha & \leqslant & \sum_{\alpha \in \N^p} \frac{C_{N}}{R^\alpha}  \frac{\prod_{i \in I} \eps_i^{N_i}}{ \prod_{j \in J} \frac{1}{\eps_j^{N_j}}} ~\rho^\alpha \\
& \leqslant & \sum_{\alpha \in \N^p} \widetilde{C}_{N}\frac{\rho^\alpha}{R^\alpha} < +\infty
\end{eqnarray*}
\end{preuve}

Thanks to this lemma, we can expand $g$:
$$g=\sum_{\alpha \in \N^p} u_{\alpha}.(T - T^0)^\alpha$$
where $u_\alpha=(u_\alpha^1,\dots,u_\alpha^{m_0})\in (\mathcal{A}_{(\overline{\C^*})^p,(\theta_I^0,\theta_J^0,t^0)}^{<0})^{m_0}$.
However, the morphism $\Psi_0$ utilizing only differential operators independent of the variable $T$, we know that $\Psi_0 (g)=0$ implies $\Psi_0(u_\alpha) =0$ for all $\alpha \in \N^p$. But, by assumption, the vector space 
$(\ker \Psi_0)_{(\theta_I^0,\theta_J^0,t^0)}$ is finite dimensional. In writing $f_1,\dots,f_d$ a basis of the latter, there exists $\lambda_1^\alpha,\dots,\lambda_d^\alpha \in \C$ such that
$$u_\alpha =\lambda_1^\alpha f_1 +\dots + \lambda_d^\alpha f_d$$ 
Then, we get $g = \sum_{\alpha \in \N^p} u_{\alpha}.(T - T^0)^\alpha = \sum_{\alpha \in \N^p}  \sum_{k=1}^d \lambda_k^\alpha f_k (T - T^0)^\alpha  = \sum_{k=1}^d \Big( f_k \sum_{\alpha \in \N^p} \lambda_k^\alpha (T - T^0)^\alpha \Big)$.\\
It remains to be shown that the series $\sum_{\alpha \in \N^p} \lambda_k^ \alpha (T - T^0)^\alpha$ converge. For that, let us show thanks to the following lemma we know to control the $ \vert \lambda_k^\alpha \vert$ using $\Vert u_{\alpha} \Vert_\infty^{K,\eps}$: 

\begin{lemme}\label{majoration-lambda-k}
There exists $C > 0$ such that for all $k=1,\dots,d$ and all $\alpha\in\N^p$, we have $\vert \lambda^\alpha_k \vert \leqslant C \Vert u_\alpha \Vert_\infty^{K,\eps}$. 
\end{lemme}

\begin{preuve}
Let $\Vert - \Vert_{\infty,f}$ be the norm on $E=\mathrm{Vect}(f_1,\dots,f_d)$ defined as follows: if $g =a_1 f_1 + \dots + a_d f_d \in E$ where $a_1,\dots,a_d \in \C$, we note
$$\Vert g \Vert_{\infty,f} := \sup_{k=1,\dots,d} \vert a_k \vert$$
Then $\vert \lambda^\alpha_k \vert \leqslant \Vert u_\alpha \Vert_{\infty,f}\leqslant C \Vert u_\alpha \Vert_\infty^{K,\eps}$ for all $k=1,\dots,d$, because of the finite dimension of $E$.
\end{preuve}

Thanks to this lemma, for $\vert T_j - T^0_j \vert = \rho_j <R_j$, we get
$$\Big\vert \sum_{\alpha \in \N^p} \lambda_k^\alpha (T - T^0)^\alpha \Big\vert \leqslant  \sum_{\alpha \in \N^p} \vert \lambda_k^\alpha \vert \rho^\alpha \leqslant C \sum_{\alpha \in \N^p} \Vert u_\alpha \Vert_\infty^{K,\eps}~\rho^\alpha$$
The convergence of this last term is guaranteed by lemma \ref{developpement}.\\

\textbf{Let us prove that $\Phi_l$ induces an isomorphism on cohomology}\\

Let us denote by $[g]$ the class of an element of the target, and $[[g]]$ the class of an element of the source.\\
Let $f_1,\dots,f_d \in (\mathcal{A}_{(\overline{\C^*})^p,(\theta_I^0,\theta_J^0,t^0)}^{<0})^{m_l}$ such that $[f_1],\dots,[f_d]$ is a basis of the vector space $\textbf{R}^l\mathcal{H}om_{\mathrm{D}}(\mathrm{M},
\mathcal{A}_{(\overline{\C^*})^p}^{<0})_{(\theta_I^0,\theta_J^0,t^0)}$.\\ 
Let $g = (g_1,\dots,g_{m_l}) \in (\mathcal{A}_{(\overline{\C^*})^p \times \C^p/\Z^p,x^0}^{<0})^{m_l}$.
Thanks to lemma \ref{developpement}, let us expand $g$ into
$$g=\sum_{\alpha \in \N^p} u_{\alpha}.(T - T^0)^\alpha$$
where $u_\alpha=(u_\alpha^1,\dots,u_\alpha^{m_l}) \in
(\mathcal{A}_{(\overline{\C^*})^p,(\theta_I^0,\theta_J^0,t^0)}^{<0})^{m_l}$.
Then, let us write:
\begin{eqnarray*}
g & = & \sum_{\alpha \in \N^p} \Big( \Psi_{l-1}(v_\alpha) + \sum_{k=1}^d \lambda_k^\alpha f_k \Big) (T - T^0)^\alpha \\
  & = &  \sum_{\alpha \in \N^p} \Psi_{l-1}(v_\alpha) (T - T^0)^\alpha + \sum_{k=1}^d f_k \sum_{\alpha \in \N^p}\lambda_k^\alpha (T - T^0)^\alpha
\end{eqnarray*}
where $v_\alpha \in (\mathcal{A}_{(\overline{\C^*})^p \times \C^p/\Z^p,x^0}^{<0})^{m_{l-1}}$. Thus, we get $\sum_{\alpha \in \N^p} \Psi_{l-1}(v_\alpha) (T - T^0)^\alpha \in \textrm{Im}(\Psi_{l-1})$.\\
Then, the morphism $\Phi_l$ induces a surjection on the cohomology groups. Indeed, we have 
\begin{eqnarray*}
[g] & = & \Big[ \sum_{\alpha \in \N^p} \Psi_{l-1}(v_\alpha) (T - T^0)^\alpha  \Big] + \Big[ \sum_{k=1}^d f_k \sum_{\alpha \in \N^p}\lambda_k^\alpha (T - T^0)^\alpha \Big] \\
     & = & [0] + \sum_{k=1}^d \Big[ f_k \sum_{\alpha \in \N^p}\lambda_k^\alpha (T - T^0)^\alpha \Big] \\
     & = & \Phi_l \Big( \sum_{k=1}^d [[f_k]] \sum_{\alpha \in \N^p}\lambda_k^\alpha (T - T^0)^\alpha  \Big)
\end{eqnarray*} 
Now, let us prove that $\Phi_l$ induces an injection on the cohomology groups.
Let us assume that $\Phi_l \Big( \sum_{k=1}^d (\sum_{\alpha \in \N^p}\lambda_k^\alpha (T - T^0)^\alpha) [[f_k]] \Big) = [0]$, namely $\Big[ \sum_{\alpha \in \N^p} (\sum_{k=1}^d \lambda_k^\alpha f_k) (T - T^0)^\alpha \Big] = [0]$, namely
$\sum_{\alpha \in \N^p} (\sum_{k=1}^d \lambda_k^\alpha f_k) (T - T^0)^\alpha \in \textrm{Im}\Psi_{l-1}$, thus $\sum_{k=1}^d \lambda_k^\alpha f_k \in \textrm{Im}\Psi_{l-1}$ for all $\alpha \in \N^p$. We get 
$[0] = \Big[ \sum_{k=1}^d \lambda_k^\alpha f_k \Big] = \sum_{k=1}^d \lambda_k^\alpha [f_k]$ then $\lambda_k^\alpha=0$ for all $\alpha \in \N^p$ because the $[f_k]$ are linearly independent on $\C$. Then, we get $\sum_{k=1}^d (\sum_{\alpha \in \N^p}\lambda_k^\alpha (T - T^0)^\alpha) [[f_k]] = [[0]]$.

					\section{A result on asymptotic expansions}

The aim of this section is to prove the following two results (propositions \ref{lemme-des-DA-1var} and \ref{lemme-des-DA-1var-parametres}). In order to achieve that, we consider the locally integrable function $K : t \flto \dfrac{-1}{2i\pi (t-1)}$ on $\C^*$ which we regard as distribution relative to the Haar measure of the multiplicative group $\C^*$. We denote by $K*f$ the partial convolution product relative to the variable $t$, namely
$$(K*f)(t,s)=\dfrac{1}{2i\pi} \int_{\C^*} f(\xi,s) \dfrac{1}{1-\frac{\xi}{t}} \dfrac{d\xi}{\xi} \wedge \dfrac{d\bar\xi}{\bar\xi}$$

\begin{prop}\label{lemme-des-DA-1var}
There exists an exact sequence of $\mathrm{D}$-modules
$$\xymatrix@R=0pt{0 \ar[r] & {(q_2 \circ \pi)}_*\mathcal{C}_{\mathbb{P}^1 \times \C}^{<0} \ar[r]^-{\partial_{\bar{t}}} & {(q_2 \circ \pi)}_*\mathcal{C}_{\mathbb{P}^1 \times \C}^{<0} \ar[r]^-{\eps} & 
\dfrac{\pi_* \mathcal{O}_{\C}^{an}[[t]][\frac{1}{t}]}{\pi_* \mathcal{O}_{\C}^{an}[\frac{1}{t}]} \oplus \dfrac{\pi_* \mathcal{O}_{\C}^{an}[[\frac{1}{t}]][t]}{(t)\pi_* \mathcal{O}_{\C}^{an}[t]} \ar[r] & 0}$$
where the image of $f$ by $\eps$ is the asymptotic expansions of $K*f$ at $0$ and infinity, and where $t$ and $t\partial_t$ act by $t$ and $t\frac{\partial}{\partial t}-s-1$.\\
Moreover: $\eps \circ (\tau t^{-1} -1) = (\tau t^{-1} -1) \circ \eps$.
\end{prop}

Let us state an analogue of proposition \ref{lemme-des-DA-1var} with parameters. This proposition will be useful in the proof of the main theorem: 

\begin{prop}\label{lemme-des-DA-1var-parametres}
For any subset $I$ of $\{1,\dots,p\}$, let us denote 
$$\mathcal{C}^{<0}_I := (q_2 \circ \pi)_*\mathcal{C}_{(\mathbb{P}^1)^{\{1,\dots,p\}\setminus I} \times \C^{p}}^{<0}$$
For all disjoint subsets $I,J \subset \{1,\dots,p\}$, the operator $\partial_{\bar{t_j}} (j\notin I \cup J)$ on
$$\dfrac{ \mathcal{C}^{<0}_{I\cup J} \Big[\Big[t_I,\dfrac{1}{t_J}\Big]\Big]\Big[\dfrac{1}{t_I},t_J\Big]}{
\Big(\displaystyle\sum_{l\in J} \big( t_l \big)  \mathcal{C}^{<0}_{I\cup J} \Big[\Big[t_I\Big]\Big]\Big[\dfrac{1}{t_I},t_J\Big] \Big) + \mathcal{C}^{<0}_{I\cup J} \Big[\Big[\dfrac{1}{t_J}\Big]\Big]\Big[t_I,t_J\Big] }$$
is injective and its cokernel is the direct sum of
$$\dfrac{ \mathcal{C}^{<0}_{I\cup J\cup\{j\}} \Big[\Big[t_j,t_I,\dfrac{1}{t_J}\Big]\Big]\Big[\dfrac{1}{t_j},\dfrac{1}{t_I},t_J\Big]}{
\Big(\displaystyle\sum_{l\in J} \big( t_l \big)  \mathcal{C}^{<0}_{I\cup J\cup\{j\}} \Big[\Big[t_j,t_I\Big]\Big]\Big[\dfrac{1}{t_j},\dfrac{1}{t_I},t_J\Big] \Big) + \mathcal{C}^{<0}_{I\cup J\cup\{j\}} \Big[\Big[\dfrac{1}{t_J}\Big]\Big]\Big[t_j,t_I,t_J\Big] }$$
and
$$\dfrac{ \mathcal{C}^{<0}_{I\cup J\cup\{j\}} \Big[\Big[t_I,\dfrac{1}{t_J},\dfrac{1}{t_j}\Big]\Big]\Big[\dfrac{1}{t_I},t_J,t_j\Big]}{
\Big(\displaystyle\sum_{l\in J\cup\{j\}} \big( t_l \big)  \mathcal{C}^{<0}_{I\cup J\cup\{j\}} \Big[\Big[t_I\Big]\Big]\Big[\dfrac{1}{t_I},t_J,t_j\Big] \Big) + \mathcal{C}^{<0}_{I\cup J\cup\{j\}}\Big[\Big[\dfrac{1}{t_J},\dfrac{1}{t_j}\Big]\Big]\Big[t_I,t_J,t_j\Big] }$$
If $\eps$ is the morphism into the cokernel of this action, the image of $f$ by $\eps$ is the asymptotic expansions of $K*f$ at $0$ and infinity relative to the variable $t_j~(j\notin I \cup J)$, and where $t_k$ and $t_k\partial_{t_k}$ act by $t_k$ and $t_k\frac{\partial}{\partial t_k}-s_k-1$.\\
Moreover: $\eps \circ (\tau_k t_k^{-1} -1) = (\tau_k t_k^{-1} -1) \circ \eps$.
\end{prop}

The proof of this proposition is similar to that of \ref{lemme-des-DA-1var}. The only point to be checked is the rapid decay condition of the coefficients of the asymptotic expansions, which we explicitly compute in the lemma \ref{DA-existence-formule}. Thus, this checking is immediate.\\
Let us prove the proposition \ref{lemme-des-DA-1var}.

\begin{lemme}\label{DA-existence-formule}
For any section $f$ of ${(q_2 \circ \pi)}_*\mathcal{C}_{\mathbb{P}^1 \times \C}^{<0}$, the convolution product $K*f$ exists and has asymptotic expansions at $0$ and at infinity which are respectively 
$$\sum_{k=1}^{+\infty} \Bigg( \frac{-1}{2i\pi} \int_{\C^*} \frac{f(\xi,s)}{\xi^k} \frac{d\xi}{\xi} \wedge \frac{d\bar\xi}{\bar\xi}\Bigg) t^k~~\textrm{et}~~\sum_{k=0}^{+\infty} \Bigg( \frac{1}{2i\pi} \int_{\C^*} \xi^k f(\xi,s) \frac{d\xi}{\xi} \wedge \frac{d\bar\xi}{\bar\xi}\Bigg) \frac{1}{t^k}$$
Moreover, $\eps$ is a morphism of $\mathrm{D}$-modules, \emph{i.e.} for any section $f$ of 
${(q_2 \circ \pi)}_*\mathcal{C}_{\mathbb{P}^1 \times \C}^{<0}$, we have $\eps\Big((t\dfrac{\partial}{\partial t}-s-1)f\Big) = (t\dfrac{\partial}{\partial t}-s-1)\eps(f)$.
\end{lemme}

\begin{preuve}
For any section $f$ of ${(q_2 \circ \pi)}_*\mathcal{C}_{\mathbb{P}^1 \times \C}^{<0}$, the convolution product 
$K*f$ is well defined because of the rapid decay condition of $f$.\\
Let us prove that $K*f$ has an asymptotic expansion at infinity. In writing 
$\displaystyle\dfrac{1}{1-\dfrac{\xi}{t}} = \sum_{k=0}^{n} \dfrac{\xi^k}{t^k} +  \dfrac{\dfrac{\xi^{n+1}}{t^{n+1}}}{1-\dfrac{\xi}{t}}$, the convolution product $K*f$ is
$$\sum_{k=0}^{n} \Big( \frac{1}{2i\pi} \int_{\C^*} \xi^k f(\xi,s) \frac{d\xi}{\xi} \wedge \frac{d\bar\xi}{\bar\xi}\Big) \frac{1}{t^k} + \Big( \frac{1}{2i\pi} \int_{\C^*}\frac{\xi^{n+1} f(\xi,s)}{1-\frac{\xi}{t}} \frac{d\xi}{\xi} \wedge \frac{d\bar\xi}{\bar\xi} \Big) \frac{1}{t^{n+1}}$$
The integrals $\displaystyle\int_{\C^*} \xi^k f(\xi,s) \frac{d\xi}{\xi} \wedge \frac{d\bar\xi}{\bar\xi}$ are finite.\\
It is the same for $\displaystyle\int_{\C^*}\frac{\xi^{n+1} f(\xi,s)}{1-\frac{\xi}{t}} \frac{d\xi}{\xi} \wedge \frac{d\bar\xi}{\bar\xi}$: in cuting out this integral as follows
$$\int_{\C^* \cap \{|\xi-t|\leqslant |t|\}} \frac{\xi^{n+1} f(\xi,s)}{1-\frac{\xi}{t}} \frac{d\xi}{\xi} \wedge \frac{d\bar\xi}{\bar\xi} + \int_{\C^* \cap \{|\xi-t| > |t|\}} \frac{\xi^{n+1} f(\xi,s)}{1-\frac{\xi}{t}} \frac{d\xi}{\xi} \wedge \frac{d\bar\xi}{\bar\xi}$$
there exists $C$ independent from $t$ such that $\Big\vert\displaystyle\int_{\C^*}\frac{\xi^{n+1} f(\xi,s)}{1-\frac{\xi}{t}} \frac{d\xi}{\xi} \wedge \frac{d\bar\xi}{\bar\xi}\Big\vert \leqslant C$ when $|t|$ is sufficiently large and $s$ in a compact.
Thus, the asymptotic expansion of $K*f$ at infinity is
$$\sum_{k=0}^{+\infty} \Big( \frac{1}{2i\pi} \int_{\C^*} \xi^k f(\xi,s) \frac{d\xi}{\xi} \wedge \frac{d\bar\xi}{\bar\xi}\Big) \frac{1}{t^k}$$
Likewise, we get the asymptotic expansion of $K*f$ at $0$.\\
Now, let us denote by $DA_\infty(g)$ and $DA_0(g)$ the asymptotic expansions of a section $g$ at infinity and at $0$, and  let us prove that $\eps\Big((t\frac{\partial}{\partial t}-s-1)f\Big) = (t\frac{\partial}{\partial t}-s-1)\eps(f)$.\\
We have $\displaystyle DA_\infty(K*t\frac{\partial f}{\partial t})=\sum_{k=0}^{+\infty} \Big( \frac{1}{2i\pi} \int_{\C^*} \xi^{k+1} \frac{\partial f}{\partial \xi}(\xi,s) \frac{d\xi}{\xi} \wedge \frac{d\bar\xi}{\bar\xi}\Big) \frac{1}{t^k}$.
In writing $\displaystyle\int_{\C^*} \xi^{k+1} \frac{\partial f}{\partial \xi}(\xi,s) \frac{d\xi}{\xi} \wedge \frac{d\bar\xi}{\bar\xi} = \int_{\C^*} \frac{\partial}{\partial \xi}(\frac{\xi^{k} f(\xi,s)}{\bar\xi}) d\xi \wedge d\bar\xi - \displaystyle\int_{\C^*} k \frac{\xi^{k-1} f(\xi,s)}{\bar\xi} d\xi \wedge d\bar\xi$ and using the Stokes theorem to show that $\displaystyle\int_{\C^*} \frac{\partial}{\partial \xi}(\frac{\xi^{k} f(\xi,s)}{\bar\xi}) d\xi \wedge d\bar\xi = 0$, we get $\displaystyle\int_{\C^*} \xi^{k+1} \frac{\partial f}{\partial \xi}(\xi,s) \frac{d\xi}{\xi} \wedge \frac{d\bar\xi}{\bar\xi} = -k \int_{\C^*} \xi^{k} f(\xi,s) \frac{d\xi}{\xi} \wedge \frac{d\bar\xi}{\bar\xi}$ for $k \geqslant 0$, thus
$$\begin{array}{ccl}
\displaystyle DA_\infty(K*t\frac{\partial f}{\partial t}) &=& \displaystyle\sum_{k=0}^{+\infty} \Big[-k\Big( \frac{1}{2i\pi} \displaystyle\int_{\C^*} \xi^{k} f(\xi,s) \frac{d\xi}{\xi} \wedge \frac{d\bar\xi}{\bar\xi} \Big)\Big] \frac{1}{t^{k}} \\
& = & t\dfrac{\partial}{\partial t} DA_\infty (K*f)
\end{array}$$
and finally $DA_\infty\Big(K*(t\dfrac{\partial}{\partial t}-s-1)f\Big)  = (t\dfrac{\partial}{\partial t}-s-1) DA_\infty (K*f)$.\\
Likewise, we get $DA_0\Big(K*(t\dfrac{\partial}{\partial t}-s-1)f\Big) = (t\dfrac{\partial}{\partial t}-s-1) DA_0 (K*f)$.\\
Thus $\eps\Big((t\dfrac{\partial}{\partial t}-s-1)f\Big) = (t\dfrac{\partial}{\partial t}-s-1)\eps(f)$.
\end{preuve}

Now, we can prove the proposition \ref{lemme-des-DA-1var}.\\

\emph{First step}: 
thanks to lemme \ref{DA-existence-formule}, the image of $\eps$ in contained in $\dfrac{\pi_* \mathcal{O}_{\C}^{an}[[t]][\frac{1}{t}]}{\pi_* \mathcal{O}_{\C}^{an}[\frac{1}{t}]} \oplus \dfrac{\pi_* \mathcal{O}_{\C}^{an}[[\frac{1}{t}]][t]}{(t)\pi_* \mathcal{O}_{\C}^{an}[t]}$.
Moreover, we get easily $\eps \Big( (\tau t^{-1} -1) f \Big) = (\tau t^{-1} -1) \eps(f)$ in the quotients above.\\

\emph{Second step}: 
If we denote by $\delta_1$ the Dirac distribution at $t=1$, we have $\dfrac{\partial K}{\partial \bar{t}}=\delta_1$. Thus, for any section $f$ of ${(q_2 \circ \pi)}_*\mathcal{C}_{\mathbb{P}^1 \times \C}^{<0}$, we get $f = \dfrac{\partial K}{\partial \bar{t}}*f = K* \dfrac{\partial f}{\partial \bar{t}}$, and the morphism $\xymatrix{{(q_2 \circ \pi)}_*\mathcal{C}_{\mathbb{P}^1 \times \C}^{<0} \ar[r]^-{\partial_{\bar{t}}} & {(q_2 \circ \pi)}_*\mathcal{C}_{\mathbb{P}^1 \times \C}^{<0}}$ is injective.\\
Moreover, we have $f = \dfrac{\partial}{\partial \bar{t}} (K*f)$. Thus, thanks to preceeding injectivity, $f$ is an element of the image of this morphism if and only if $K*f$ satisfy the rapid decay condition at $0$ and at infinity, \emph{i.e.} $f \in Ker(\eps)$.\\

\emph{Third step}:
It remains to be shown that $\eps$ is surjective. Let $\gamma_0$ and $\gamma_\infty$ be sections of $\dfrac{\pi_* \mathcal{O}_{\C}^{an}[[t]][\frac{1}{t}]}{\pi_* \mathcal{O}_{\C}^{an}[\frac{1}{t}]} \oplus \dfrac{\pi_* \mathcal{O}_{\C}^{an}[[\frac{1}{t}]][t]}{(t)\pi_* \mathcal{O}_{\C}^{an}[t]}$. thanks to Borel theorem (see \cite{Mal-livre} p.62 or \cite{Mal-pt-sing-irreg} or \cite{Ram}), there exists two functions $g_0$ and $g_\infty$ which are $\mathcal{C}^{\infty}$ in $t$ on $\mathbb{P}^1$ (and holomorphic locally in $s$) and which have respectively $\gamma_0$ for asymptotic expansion at $0$, and $\ gamma_\infty$ for asymptotic expansion at infinity. Then there exists a function $g$ which is $\mathcal{C}^{\infty}$ in $t$ on $\mathbb{P}^1$ (and holomorphic locally in $s$) and such that $\gamma_0$ and $\gamma_\infty$ are the asymptotic expansions of $g$ at $0$ and infinity.\\
The section $\dfrac{\partial g}{\partial \bar{t}}$ is a section of ${(q_2 \circ \pi)}_*\mathcal{C}_{\mathbb{P}^1 \times \C}^{<0}$
because $\bar{t}$ don't appear in the asymptotic expansions of $g$. Moreover, $K*\dfrac{\partial g}{\partial \bar{t}} =\dfrac{\partial K}{\partial \bar{t}}*g=g$. Thus, the asymptotic expansions of $K*\dfrac{\partial g}{\partial \bar{t}}$ are those of $g$, \emph{i.e.} $\gamma_0$ and $\gamma_\infty$, and $\dfrac{\partial g}{\partial \bar{t}}$ is an antecedent of $(\gamma_0,\gamma_\infty)$ by $\eps$.

					\section{The proof of the main theorem}

In this section, we prove the theorem \ref{th-mellin}.

\subsection{Construction of the kernel $\noyau$}

We denote by $\widetilde{\mathrm{D}}$ the non-commutative $\C$-algebra $\mathrm{D}[s_1, \dots ,s_p]\langle \tau_1,{\tau_1}^{-1}, \dots,\tau_p ,\tau_p^{-1}\rangle$ generated by the $t_j$, $t_j^{-1}$, $t_j\partial_{t_j}$, $s_j$, $\tau_j$, $\tau_j^{-1}$ ($j$ from $1$ to $p$) and the relations
$$(t_j\partial_{t_j})t_j - t_j(t_j\partial_{t_j})=t_j$$
$$\tau_j s_j =(s_j+1) \tau_j$$
$$s_j t_j=t_j s_j$$
$$(t_j\partial_{t_j})s_j=s_j (t_j\partial_{t_j})$$
$$\tau_j t_j = t_j \tau_j$$
$$(t_j\partial_{t_j}) \tau_j=\tau_j (t_j\partial_{t_j})$$

First of all, let us define a sheaf $\noyau$ which plays the role of the Mellin transform kernel.
In order to do that, let us consider the natural $\widetilde{\mathrm{D}}$-module structure on the sheaf $\overline{\mathcal{O}}=k_*\pi_*\mathcal{O}_{(\C^*)^p \times \C^p}^{an}$.

\begin{defi}
Let $\mathbb{T}^{s+1}$ be the $\widetilde{\mathrm{D}}$-module
$$\dfrac{\widetilde{\mathrm{D}}}{\widetilde{\mathrm{D}}\Big(t_1\frac{\partial}{\partial t_1}-s_1 -1,\dots,t_p\frac{\partial}{\partial t_p}-s_p -1,\tau_1 t^{-1}_1 -1,\dots,\tau_p t^{-1}_p -1\Big)}$$
The sheaf $\noyau$ is defined as
\begin{center}
$\noyau=\mathcal{H}om_{\widetilde{\mathrm{D}}} \Big(  \mathbb{T}^{s+1} ,
\overline{\mathcal{O}} \Big)=\mathcal{H}om_{\widetilde{\mathrm{D}}} \Big(  \mathbb{T}^{s+1} ,
k_*\pi_*\mathcal{O}_{(\C^*)^p \times \C^p}^{an} \Big)$
\end{center}
where $\widetilde{\mathrm{D}}$ and $\mathbb{T}^{s+1}$ are constant sheaves on $(\overline{\C^*})^p \times \C^p/\Z^p$.\\
\end{defi}

Locally on $(\overline{\C^*})^p \times \C^p/\Z^p$, we can choose a determination of $\tsp=t_1^{s_1 +1}\cdots t_p^{s_p +1}$. A local section of $\noyau$ consists of the data of a local section $\varphi$ of $\overline{\mathcal{O}}$ which satisfies the equations defining $\mathbb{T}^{s+1}$. Consequently, $\noyau$ is the $r_{2}^{-1}\mathcal{O}_{\C^p / \Z^p}$-module $r_{2}^{-1}\mathcal{O}_{\C^p / \Z^p}.\tsp \subset \overline{\mathcal{O}}$. In particular, we get:

\begin{propriété}
Following the terminology of Deligne (\cite{Del}), the sheaf $\noyau$ on $(\overline{\C^*})^p \times \C^p/\Z^p$ is a \emph{relative local system} of rank $1$ free $r_{2}^{-1} \mathcal{O}_{\C^p / \Z^p}$-modules, \emph{i.e.} on any open set $U \times V$ sufficiently small, there exists an isomorphism $\noyau \vert_{U\times V} \simeq r_{2}^{-1}\mathcal{O}_V$. 
\end{propriété}

Following results of Deligne (\cite{Del}), there is an equivalence between the category of relative local systems of $r_{2}^{-1} \mathcal{O}_{\C^p / \Z^p}$-modules on $(\overline{\C^*})^p \times \C^p/\Z^p$ and that of $\mathcal{O}_{\C^p / \Z^p}\big[\pi_1\big((\overline{\C^*})^p\big)\big]$-modules. If we denote $\C\big[\pi_1\big((\overline{\C^*})^p\big)\big]=\C[M,M^{-1}]$ where $M=(M_1,\dots,M_p)$, we get $\mathcal{O}_{\C^p / \Z^p}\big[\pi_1\big((\overline{\C^*})^p\big)\big]=\mathcal{O}_{\C^p / \Z^p}[M,M^{-1}]$.\\
Since 
$$M_j e^{\sigma_1(s_1 + 1)}\cdots e^{\sigma_p(s_p + 1)}  = T_j^{-1} e^{\sigma_1(s_1 + 1)}\cdots e^{\sigma_p(s_p + 1)}$$
the sheaf $\noyau$ corresponds to the $\mathcal{O}_{\C^p / \Z^p}[M,M^{-1}]$-module $\mathcal{O}_{\C^p / \Z^p}$ on which $M_j$ acts by the multiplication by $T_j^{- 1}$, \emph{via} the equivalence recalled above.\\

We consider the constant sheaf of rings $\C[T,T^{-1}]$ on $\C^p / \Z^p$ as a subsheaf of $\mathcal{O}_{\C^p / \Z^p}$. We still denote by $\C[T,T^{-1}]$ its inverse image by $r_2$ and $r_{2} \circ \bar{q}$.\\
Let us recall that $\bar{\mathcal{L}}$ is the natural local system of rank $1$ free $\C[M,M^{-1}]$-modules on $(\overline{\C^*})^p$ defined in section \ref{def-mellin-sheaves}, \emph{i.e.} $\bar{q}_!\bar{q}^{-1}\underline{\C}_{(\overline{\C^*})^p}$. \emph{Via} the equivalence recalled above, the sheaf $r_1^{-1}\bar{\mathcal{L}}$ equipped with the trivial structure of $r_{2}^{-1} \mathcal{O}_{\C^p / \Z^p}$-module corresponds to the $\C[M,M^{-1}]$-module $\C[T,T^{-1}]$ on which $M_j$ acts by the multiplication by $T_j^{- 1}$.\\
Likewise, $r_{2}^{-1} \mathcal{O}_{\C^p / \Z^p}$ corresponds to the $\mathcal{O}_{\C^p / \Z^p}[M,M^{-1}]$-module $\mathcal{O}_{\C^p / \Z^p}$ on which $M_j$ acts by $1$.\\
We deduce from it that the $r_{2}^{-1} \mathcal{O}_{\C^p / \Z^p}$-module $r_1^{-1}\bar{\mathcal{L}} \otimes_{\C[T,T^{-1}]} r_{2}^{-1} \mathcal{O}_{\C^p / \Z^p}$ corresponds to the $\mathcal{O}_{\C^p / \Z^p}[M,M^{-1}]$-module $\C[T,T^{-1}] \otimes_{\C[T,T^{-1}]} \mathcal{O}_{\C^p / \Z^p}$ on which $M_j$ acts by $M_j.(a\otimes b) = T_j^{-1}a\otimes 1b = T_j^{-1}a\otimes b = a\otimes T_j^{-1}b$.\\

Now, let us note this isomorphism in the category of $\mathcal{O}_{\C^p / \Z^p}[M,M^{-1}]$-modules:
$$\begin{array}{rcl}
\C[T,T^{-1}] \otimes_{\C[T,T^{-1}]} \mathcal{O}_{\C^p / \Z^p} & \stackrel{\sim}{\fl} & \mathcal{O}_{\C^p / \Z^p} ~~~~~~(\sharp)\\
1 \otimes f & \flto & f
\end{array}$$ 
where, $M_j$ acts by the multiplication by $T_j^{- 1}$ on $\mathcal{O}_{\C^p / \Z^p}$ and by the multiplication by $T_j^{-1}$ on the factor $\C[T,T^{-1}]$ of the term $\C[T,T^{-1}] \otimes_{\C[T,T^{-1}]} \mathcal{O}_{\C^p / \Z^p}$.\\
\emph{Via} the equivalence recalled above, this isomorphism corresponds to an isomorphism in the category of relative local systems of $r_{2}^{-1}\mathcal{O}_{\C^p / \Z^p}$-modules 
$$r_1^{-1}\bar{\mathcal{L}} \otimes_{\C[T,T^{-1}]} r_{2}^{-1}\mathcal{O}_{\C^p / \Z^p} \stackrel{\sim}{\fl} \noyau~~~~~~(\sharp\sharp)$$
or 
$$\dfrac{r_1^{-1}\bar{\mathcal{L}} \otimes_{\C} r_{2}^{-1}\mathcal{O}_{\C^p / \Z^p}}{(M_1T_1 - 1, \dots ,M_p T_p - 1 )} \stackrel{\sim}{\fl} \noyau$$
However, let us note that if we compose $(\sharp)$ with the multiplication by a unit $T^\alpha~(\alpha\in\Z^p)$ of $\C[T,T^{-1}]$, we get an other isomorphism
$$\begin{array}{rcl}
\C[T,T^{-1}] \otimes_{\C[T,T^{-1}]} \mathcal{O}_{\C^p / \Z^p} & \stackrel{\sim}{\fl} & \mathcal{O}_{\C^p / \Z^p} \\
1 \otimes f & \flto & T^\alpha f
\end{array}$$
Thus, the isomorphism $(\sharp\sharp)$ is not unique.\\

\subsection{Another formulation of the Mellin transform on sheaves}

The definition of the Mellin transform of a sheaf $\F$ in \ref{def-mellin-faisceautique} is an algebraic definition: it uses the Alexander modules of $\F$ which are $\C[T,T^{-1}]$-modules. In the proof of theorem \ref{th-mellin}, we will use a more analytic formulation which translates better than the former one the classical definiton of the integral Mellin transform.

\begin{prop}\label{nouvelle-def}
For any object $\F$ in $D^b((\overline{\C^*})^p,\C)$, we have
$$\mathfrak{M}(\F) \cong \textbf{R} {r_2}_* ( r_{1}^{-1}\F \ol_{\C} \noyau )[p]$$
\end{prop}

\begin{preuve}
Let us calculate:
$$\begin{array}{cl}
& \textbf{R} {r_2}_* ( r_{1}^{-1}\F \ol_{\C} \noyau ) \\
= & \textbf{R} {r_2}_* \big( r_{1}^{-1} (\F \ol_{\C} \bar{\mathcal{L}} ) \otimes_{\C[T,T^{-1}]} r_{2}^{-1}\mathcal{O}_{\C^p / \Z^p} \big) \\
= & \textbf{R} {r_2}_* r_{1}^{-1} ( \F \ol_{\C} \bar{\mathcal{L}} ) \otimes_{\C[T,T^{-1}]} \mathcal{O}_{\C^p / \Z^p} ~\textrm{because $r_2$ is proper}\\
= & \textbf{R}\Gamma \Big( (\overline{\C^*})^p , \F \ol_{\C} \bar{\mathcal{L}} \Big) \otimes_{\C[T,T^{-1}]} \mathcal{O}_{\C^p / \Z^p} \\
= & \mathfrak{M}(\F)[-p]
\end{array}$$
\end{preuve}

\subsection{Computation of the complex $\mathfrak{M}(\mathcal{S}ol^{<0}(\mathcal{M}))$}

\begin{prop}\label{prop-mellin-sol}
For any coherent algebraic $\D_{(\C^*)^p}$-module $\mathcal{M}$, such that 
$\mathcal{S}ol^{<0}(\mathcal{M})\big\vert_{(\C^*)^p}$ has $\C$-constructible cohomology and $\mathcal{S}ol^{<0}(\mathcal{M})\big\vert_{(\overline{\C^*})^p \setminus (\C^*)^p}$ has $\R$-constructible cohomology, and whose module of global sections is denoted by $\mathrm{M}$,
there exists a canonical isomorphism in $D^b(\mathfrak{Mod}(\mathcal{O}_{\C^p / \Z^p}))$
$$\mathfrak{M}(\mathcal{S}ol^{<0}(\mathcal{M}))  \cong 
\textbf{R}\mathcal{H}om_{\mathrm{D}} \Big(  \mathrm{M} ,
\textbf{R}{r_2}_* ( \mathcal{A}_{(\overline{\C^*})^p \times \C^p/\Z^p}^{<0} \tnoyau) \Big)[p]$$ 
where $t_j\partial_{t_j}$ acts by
$t_j\partial_{t_j}.g(\underline{t},\underline{T}) \otimes \omega = t_j\dfrac{\partial g}{\partial t_j}(\underline{t},\underline{T}) \otimes \omega$ on the sheaf \\ $\textbf{R}{r_2}_* ( \mathcal{A}_{(\overline{\C^*})^p \times \C^p/\Z^p}^{<0} \tnoyau)$. 
\end{prop}

\begin{preuve}
thanks to proposition \ref{nouvelle-def} and the lemma \ref{sol-dec-rapide}, we get
$$\mathfrak{M}(\mathcal{S}ol^{<0}(\mathcal{M})) \cong \textbf{R} {r_2}_* ( r_{1}^{-1}\textbf{R}\mathcal{H}om_{\mathrm{D}}
( \mathrm{M} , \mathcal{A}_{(\overline{\C^*})^p}^{<0} ) \ol_{\C} \noyau )~[p]$$ 
On the one hand, we have a canonical isomorphism\\
${r_1}^{-1}\textbf{R}\mathcal{H}om_{\mathrm{D}} ( \mathrm{M} , \mathcal{A}_{(\overline{\C^*})^p}^{<0}) \ol_{\C} \noyau$

\hfill $\cong  
\textbf{R}\mathcal{H}om_{\mathrm{D}} ( \mathrm{M} , {r_1}^{-1}\mathcal{A}_{(\overline{\C^*})^p}^{<0}) \ol_{\C} r_{2}^{-1}\mathcal{O}_{\C^p / \Z^p} \ol_{r_{2}^{-1}\mathcal{O}_{\C^p / \Z^p}} \noyau$\\
and, thanks to theorem \ref{finitude-solutions}, we get a canonical isomorphism
$${r_1}^{-1}\textbf{R}\mathcal{H}om_{\mathrm{D}} ( \mathrm{M} , \mathcal{A}_{(\overline{\C^*})^p}^{<0}) \ol_{\C} \noyau \cong  \textbf{R}\mathcal{H}om_{\mathrm{D}} ( \mathrm{M} , \mathcal{A}_{(\overline{\C^*})^p \times \C^p/\Z^p}^{<0}) \ol_{r_{2}^{-1}\mathcal{O}_{\C^p / \Z^p}} \noyau$$
On the other hand, we have a canonical isomorphism\\
$\textbf{R}\mathcal{H}om_{\mathrm{D}} ( \mathrm{M} ,
\mathcal{A}_{(\overline{\C^*})^p \times \C^p/\Z^p}^{<0}) \ol_{r_{2}^{-1}\mathcal{O}_{\C^p / \Z^p}} \noyau$

\hfill $\cong \textbf{R}\mathcal{H}om_{\mathrm{D}} ( \mathrm{M} ,
\mathcal{A}_{(\overline{\C^*})^p \times \C^p/\Z^p}^{<0} \otimes_{r_{2}^{-1}\mathcal{O}_{\C^p / \Z^p}} \noyau)$\\
where the action of $\mathrm{D}$ on $\mathcal{A}_{(\overline{\C^*})^p \times \C^p/\Z^p}^{<0} \otimes_{r_{2}^{-1}\mathcal{O}_{\C^p / \Z^p}} \noyau$ is the action only on $\mathcal{A}_{(\overline{\C^*})^p \times \C^p/\Z^p}^{<0}$.\\
Finally, we get the wished isomorphism by applying the functor $\textbf{R}{r_2}_*(-)[p]$ and noting that\\

$\textbf{R}{r_2}_* \textbf{R}\mathcal{H}om_{\mathrm{D}} \Big(  \mathrm{M} , \mathcal{A}_{(\overline{\C^*})^p \times \C^p/\Z^p}^{<0} \tnoyau \Big)$

\hfill $\cong 
\textbf{R}\mathcal{H}om_{\mathrm{D}} \Big(  \mathrm{M} , \textbf{R}{r_2}_* ( \mathcal{A}_{(\overline{\C^*})^p \times \C^p/\Z^p}^{<0} \tnoyau) \Big)$\\

because $\mathrm{D}$ and $\mathrm{M}$ are constant sheaves.
\end{preuve}

However, we can express the $\mathrm{D}$-module $\textbf{R}{r_2}_* ( \mathcal{A}_{(\overline{\C^*})^p \times \C^p/\Z^p}^{<0} \tnoyau)$ where $t_j\partial_{t_j}$ acts by
$t_j\partial_{t_j}.g(\underline{t},\underline{T}) \otimes \omega = t_j\dfrac{\partial g}{\partial t_j}(\underline{t},\underline{T}) \otimes \omega$ in another way.

\subsubsection{Solving a system of difference equations}

Let us consider the injection
$$\xymatrix@R=0pt{
\mathcal{A}_{(\overline{\C^*})^p  \times \C^p/\Z^p}^{<0} \tnoyau \ar@{^{(}->}[r]^-{\iota} & \pi_*\mathcal{A}_{(\overline{\C^*})^p  \times \C^p}^{<0} \\
g(t_1,\dots,t_p,T_1,\dots,T_p) \otimes \omega \ar@{|->}[r] & g\big(t_1,\dots,t_p,e^{-2i\pi s_1},\dots,e^{-2i\pi s_p}\big) \omega}$$
where all sections are viewed as sections of $k_*\pi_*\mathcal{O}_{(\C^*)^p \times \C^p}^{an}$.\\
The injection $\iota$ is a morphism of $\mathrm{D}$-modules where $t_j\partial_{t_j}$ acts on $\pi_*\mathcal{A}_{(\overline{\C^*})^p \times \C^p}^{<0}$ by the derivation $t_j\dfrac{\partial}{\partial{t_j}} -s_j -1$ and the struture of $\mathrm{D}$-module of $\mathcal{A}_{(\overline{\C^*})^p  \times \C^p/\Z^p}^{<0} \tnoyau$ is that given to the proposition \ref{prop-mellin-sol}. Indeed, by definition of $\noyau$, we know that $t_j\dfrac{\partial \omega}{\partial t_j}=(s_j+1)\omega$ for any section $\omega$ of $\noyau$, and:
$$\begin{array}{rl} 
  & (t_j\dfrac{\partial}{\partial{t_j}} -s_j -1)\iota\big( g(\underline{t},\underline{T}) \otimes \omega \big) \\
=& (t_j\dfrac{\partial}{\partial{t_j}} -s_j -1).g(\underline{t},\underline{T})\omega \\
=& t_j\dfrac{\partial g}{\partial t_j}(\underline{t},\underline{T})\omega + (s_j +1)g(\underline{t},\underline{T})\omega - (s_j +1)g(\underline{t},\underline{T})\omega  \\
=& \iota \big( t_j\dfrac{\partial g}{\partial t_j}(\underline{t},\underline{T}) \otimes \omega) \big)  \\
=& \iota \big( t_j\partial_{t_j}~.~g(\underline{t},\underline{T}) \otimes \omega) \big)
\end{array}$$

Let $\mathcal{K}^\bullet (\tau_1 t_1^{-1} -1, \dots ,\tau_p t_p^{-1} -1;\pi_*\mathcal{A}_{(\overline{\C^*})^p \times \C^p}^{<0})$ be the Koszul complex where the degree $0$ is placed on the left.

\begin{prop}\label{tau-t-dec-rapide-pvar}
There is an isomorphism in the category of $\mathrm{D}$-modules:
$$\mathcal{K}^\bullet (\tau_1 t_1^{-1} -1, \dots ,\tau_p t_p^{-1} -1;\pi_*\mathcal{A}_{(\overline{\C^*})^p \times \C^p}^{<0}) \cong \mathcal{A}_{(\overline{\C^*})^p  \times \C^p/\Z^p}^{<0} \tnoyau$$
where $t_j\partial_{t_j}$ acts on the Koszul complex by the derivation $t_j\dfrac{\partial}{\partial{t_j}} -s_j -1$ and the struture of $\mathrm{D}$-module of $\mathcal{A}_{(\overline{\C^*})^p  \times \C^p/\Z^p}^{<0} \tnoyau$ is that given to the proposition \ref{prop-mellin-sol}.
\end{prop}

\begin{preuve}
This Koszul complex is viewed as the total complex associated with the complex of order $p$ constituted of the $p$ complexes
$\xymatrix@R=0pt{[\pi_*\mathcal{A}_{\overline{(\C^*)}^p \times \C^p}^{<0} \ar[r]^-{\tau_j t_j^{-1} -1} & \pi_*\mathcal{A}_{\overline{(\C^*)}^p \times \C^p}^{<0}] \\
                                                                        \bullet &  }$, and we proceed by induction.\\
For that, let us establish some notations. If $I$ is a subset of $\{1,\dots,p\}$, let us denote $\mathbb{T}_I^{s+1}=\dfrac{\widetilde{\mathrm{D}}}{\widetilde{\mathrm{D}}\Big(t_{I}\frac{\partial}{\partial t_{I}}-s_{I} -1,\tau_{I} t^{-1}_{I} -1\Big)}$ and
$$\noyau_I:=\mathcal{H}om_{\widetilde{\mathrm{D}}} \Bigg(  \mathbb{T}_I^{s+1} ,
k_*\pi_*\mathcal{O}_{(\C^*)^p \times \C^p}^{an} \Bigg)$$ 
                                                               
Let us treat the initial case: the morphism 
$$\xymatrix@C=2cm{\pi_*\mathcal{A}_{(\overline{\C^*})^p   \times \C^p}^{<0} \ar[r]^-{(\tau_{1} t_{1}^{-1} -1).} & \pi_*\mathcal{A}_{(\overline{\C^*})^p  \times \C^p}^{<0} }$$
is surjective and, by definition of $\noyau_{\{1\}}$, its kernel is
$$\pi_*\mathcal{A}_{(\overline{\C^*})^p  \times \C/\Z \times \C^{\{2,\dots,p\}}}^{<0} \tnoyau_{\{1\}} $$
To establish the surjectivity, we consider the stalk at a point $(t_0,T_0)$ and, if $g$ is a section of $(\pi_*
\mathcal{A}_{\overline{\C^*} \times \C}^{<0})_{(t_0,T_0)}$, we look for $u$ in $(\pi_*\mathcal{A}_{\overline{\C^*} \times \C}^{<0})_{(t_0,T_0)}$ such that $(\tau t^{-1} -1)u=g$. In writing $v:=t^{-s-1} u$, we just have to prove the surjectivity of 
$$\xymatrix{\pi_*\mathcal{A}_{\overline{\C^*} \times \C}^{<0} \ar[r]^-{(\tau_{1} -1).} & \pi_*\mathcal{A}_{\overline{\C^*} \times \C}^{<0}}$$
rapid decay conditions being preserved through $\ts$.\\ 
We treat this last problem on small open sets.\\

Let us continue the induction.\\ 
The sheaf $\noyau_I$ is a relative local system of rank $1$ free $r_{2}^{-1} \mathcal{O}_{\C^p / \Z^p}$-modules.
It is locally generated by a determination of $t_I^{s_I+1}$.\\ 
Now, let us denote $I=\{1,\dots,j-1\}\subset\{1,\dots,p-1\}$ and $J=\{j,\dots,p\}$.\\
As above, the image of the injective natural morphism
$$\xymatrix@R=0pt{\pi_*\mathcal{A}_{(\overline{\C^*})^p  \times (\C/\Z)^{I\cup\{j\}} \times \C^{J\setminus\{j\}}}^{<0} \tnoyau_{\{j\}} \ar@{^{(}->}[r] & \pi_*\mathcal{A}_{(\overline{\C^*})^p  \times (\C/\Z)^{I} \times \C^{J}}^{<0}  \\
g(t,T_{I\cup\{j\}},s_{J\setminus\{j\}}) \otimes \omega \ar@{|->}[r] & g(t,T_I,e^{-2i\pi s_j}_j,s_{J\setminus\{j\}}) \omega  }$$
is the kernel of the surjective morphism
$$\xymatrix@C=2cm{\pi_*\mathcal{A}_{(\overline{\C^*})^p  \times (\C/\Z)^{I} \times \C^{J}}^{<0} \ar[r]^-{(\tau_j t_j^{-1} -1).} & \pi_*\mathcal{A}_{(\overline{\C^*})^p  \times (\C/\Z)^{I} \times \C^{J}}^{<0}}$$
However, $\noyau_I$ and $\noyau_{\{j\}}$ being locally $r_{2}^{-1}\mathcal{O}_{(\C / \Z)^p}$-free, we get
$$ \noyau_I \otimes_{r_{2}^{-1}\mathcal{O}_{(\C / \Z)^p}} \noyau_{\{j\}} = \noyau_{I \cup \{j\}}$$
and the functor $- \otimes_{r_{2}^{-1}\mathcal{O}_{(\C / \Z)^p}} \noyau_I$ is exact.
Thus, the image of the injective natural morphism
$$\xymatrix{\pi_*\mathcal{A}_{(\overline{\C^*})^p  \times (\C/\Z)^{I\cup\{j\}} \times \C^{J\setminus\{j\}}}^{<0} \tnoyau_{I \cup \{j\}} \ar@{^{(}->}[r] & \pi_*\mathcal{A}_{(\overline{\C^*})^p  \times (\C/\Z)^{I} \times \C^{J}}^{<0} \tnoyau_I }$$
is the kernel of the surjective morphism
$$\xymatrix@C=2cm{\pi_*\mathcal{A}_{(\overline{\C^*})^p  \times (\C/\Z)^{I} \times \C^{J}}^{<0} \tnoyau_I \ar[r]^-{(\tau_j t_j^{-1} -1).} & \pi_*\mathcal{A}_{(\overline{\C^*})^p  \times (\C/\Z)^{I} \times \C^{J}}^{<0} \tnoyau_I}$$
\end{preuve}

\subsubsection{Computation of $\textbf{R}{r_2}_*\mathcal{K}^\bullet (\tau_1 t_1^{-1} -1, \dots ,\tau_p t_p^{-1} -1;\pi_*\mathcal{A}_{(\overline{\C^*})^p \times \C^p}^{<0})$}

Let us denote 
$$\mathcal{K}=\textbf{R} {r_2}_*\mathcal{K}^\bullet ( \tau_1 t_1^{-1} -1, \dots ,\tau_p t_p^{-1} -1 ;\pi_*\mathcal{A}_{(\overline{\C^*})^p \times \C^p}^{<0})$$
with the structure of $\mathrm{D}$-module recalled above.\\                                                                                                    

By section \ref{fast-dec}, we consider the Dolbeault resolution 
$\pi_*\mathcal{A}_{(\overline{\C^*})^p \times \C^p}^{<0} \cong \pi_*\mathcal{P}^{0,\bullet}_{(\overline{\C^*})^p \times \C^p}$
which is ${r_2}_*$-acyclic. Thus, 
$$\mathcal{K} \cong \mathcal{T}ot^\bullet \Big(\mathcal{K}^\bullet (\tau_1 t_1^{-1} -1, \dots ,\tau_p t_p^{-1} -1;{(r_2 \circ \pi)}_*\mathcal{P}^{0,\bullet}_{(\overline{\C^*})^p \times \C^p})\Big)$$
In writing ${(r_2 \circ \pi)}_*={(q_2 \circ \pi)}_* \bar\pi_*$, we get 
$$\mathcal{K} \cong \mathcal{T}ot^\bullet \Big(\mathcal{K}^\bullet (\tau_1 t_1^{-1} -1, \dots ,\tau_p t_p^{-1} -1;{(q_2 \circ \pi)}_*\mathcal{P}^{0,\bullet}_{(\mathbb{P}^1)^p \times \C^p})\Big)$$

In the continuation, if $\Lambda$ is a set of indices and $R$ a ring, let $R[[t_\Lambda]]$ be the ring of formal power series with all variables $t_\lambda~(\lambda\in\Lambda)$.\\
Let us denote 
$$\mathcal{DA}_{I,J} =\dfrac{ \pi_* \mathcal{O}_{\C^p} \Big[\Big[t_I,\dfrac{1}{t_J}\Big]\Big]\Big[\dfrac{1}{t_I},t_J\Big]}{
\Big(\displaystyle\sum_{l\in J} \big( t_l \big)  \pi_* \mathcal{O}_{\C^p} \Big[\Big[t_I\Big]\Big]\Big[\dfrac{1}{t_I},t_J\Big] \Big) + \pi_* \mathcal{O}_{\C^p} \Big[\Big[\dfrac{1}{t_J}\Big]\Big]\Big[t_I,t_J\Big] }$$
for all partitions (eventually trivial) $I \coprod J = \{1, \dots, p\}$, and 
$$ \mathcal{DA} = \bigoplus_{I,J} \mathcal{DA}_{I,J}$$
where the direct sum is done with all partitions (eventually trivial) $I \coprod J = \{1, \dots, p\}$.

We get:

\begin{lemme}\label{decomposition-DA}
There exists an isomorphism
$${(q_2 \circ \pi)}_*\mathcal{P}^{0,\bullet}_{(\mathbb{P}^1)^p \times \C^p} \cong \mathcal{DA}~[-p]$$
in the category of complexes of $\mathrm{D}$-modules, where $t_j$ and $t_j\partial_{t_j}$ act by $t_j$ and $t_j\frac{\partial}{\partial t_j}-s_j-1$.
\end{lemme}

\begin{preuve}
We have ${(q_2 \circ \pi)}_*\mathcal{P}^{0,\bullet}_{(\mathbb{P}^1)^p \times \C^p} \cong \mathcal{K}^\bullet (\partial_{\bar{t_1}}, \dots , \partial_{\bar{t_p}} ; {(q_2 \circ \pi)}_*\mathcal{C}_{(\mathbb{P}^1)^p \times \C^p}^{<0})$.
This Koszul complex is viewed as the total complex associated with the complex of order $p$ constituted of the $p$ complexes
$$\xymatrix@R=0pt{[{(q_2 \circ \pi)}_* \mathcal{C}^{<0}_{(\mathbb{P}^1)^p \times \C^p} \ar[r]^-{\partial_{\bar{t_j}}} & {(q_2 \circ \pi)}_* \mathcal{C}^{<0}_{(\mathbb{P}^1)^p \times \C^p}] \\
                                                                        \bullet &  }$$ 
and we proceed by induction thanks to proposition \ref{lemme-des-DA-1var-parametres}.
\end{preuve}

Hence, we get
$$\mathcal{K} \cong \mathcal{K}^\bullet (\tau_1 t_1^{-1} -1, \dots ,\tau_p t_p^{-1} -1 ; \mathcal{DA} ) ~[-p]$$
\emph{i.e.}
$$\mathcal{K} \cong  \bigoplus_{I,J} \mathcal{K}^\bullet \Bigg(\tau_1 t_1^{-1} -1, \dots ,\tau_p t_p^{-1} -1 ; \mathcal{DA}_{I,J} \Bigg)~[-p]$$ 
in the category of complexes of $\mathrm{D}$-modules, where the direct sum is done with all partitions (eventually trivial) $I \coprod J = \{1, \dots, p\}$.

Let us simplify this complex.

\begin{lemme}\label{DAnul-pvar}
If $J \neq \emptyset$ then 
$$\mathcal{K}^\bullet \Bigg(\tau_1 t_1^{-1} -1, \dots ,\tau_p t_p^{-1} -1 ;  \mathcal{DA}_{I,J} \Bigg) \cong 0 $$
\end{lemme}

\begin{preuve}
Let us fix a nonempty set $J$.\\
This Koszul complex is viewed as the total complex associated with the complex of order $p$ constituted of the $p$ complexes 
$\xymatrix@R=0pt@C=2cm{[\mathcal{DA}_{I,J} \ar[r]^-{\tau_j t_j^{-1} -1} & \mathcal{DA}_{I,J}] \\
                                                                        \bullet &  }$ which we calculate a variable after the other.\\ 
Let $m$ be $\max(J)$.\\
\emph{Case A : } $m=p$\\
Let us denote by $\dsum_{n=0}^{\infty} \dfrac{a_n}{{t_p}^n}$ a section of $\mathcal{DA}_{I,J}$. In this quotient, we get 
$$(\tau_p {t_p}^{-1} -1)\dsum_{n=0}^{\infty} \dfrac{a_n}{{t_p}^n} =-a_0+ \dsum_{n=1}^{\infty} \dfrac{\tau_p a_{n-1}-a_{n}}{{t_p}^n}$$ 
Hence, the action of $\tau_p t_p^{-1} -1$ on $\mathcal{DA}_{I,J}$ is bijective and we get
$$\mathcal{K}^\bullet (\tau_1 t_1^{-1}, \dots ,\tau_{p-1} t_{p-1}^{-1} -1 ; 0) = 0$$ 
\emph{Case B :} $m<p$\\
Let us denote by $\dsum_{n=1}^{\infty} b_n {t_p}^n$ a section of $\mathcal{DA}_{I,J}$. In this quotient, we get
$$(\tau_p {t_p}^{-1} -1)\dsum_{n=1}^{\infty} b_n {t_p}^n = \dsum_{n=1}^{\infty} (\tau_p b_{n+1}-b_n) {t_p}^n$$
Hence, the action of $\tau_p t_p^{-1} -1$ on $\mathcal{DA}_{I,J}$ is surjective and its kernel is $\mathcal{DA}_{I-\{p\},J}$.\\
Thus, the searched complex is
$$\mathcal{K}^\bullet \Bigg(\tau_1 t_1^{-1} -1, \dots ,\tau_{p-1} t_{p-1}^{-1} -1 ; \mathcal{DA}_{I-\{p\},J} \Bigg)$$
By reiterating the process, we get
$$\mathcal{K}^\bullet \Bigg(\tau_1 t_1^{-1} -1, \dots ,\tau_m t_m^{-1} -1 ; \mathcal{DA}_{I-\{m+1, \dots ,p\},J} \Bigg)$$
which we treat as the case A because $m=\max(J)$. Hence, the searched complex is null.
\end{preuve}

However, if $J=\emptyset$, we get:

\begin{lemme}\label{DAnonnul-pvar}
There exists an isomorphism of complexes of $\mathrm{D}$-modules
$$\mathcal{K}^\bullet \Bigg(\tau_1 t_1^{-1} -1, \dots ,\tau_p t_p^{-1} -1 ; 
\dfrac{ \pi_* \mathcal{O}_{\C^p} \Big[\Big[t_1,\dots,t_p\Big]\Big]\Big[\dfrac{1}{t_1},\dots,\dfrac{1}{t_p}\Big]}{
\pi_* \mathcal{O}_{\C^p} \Big[t_1,\dots,t_p\Big] }
 \Bigg) \cong \pi_*
\mathcal{O}_{\C^p}$$ 
where $t_j$ and $t_j\partial_{t_j}$ act on the left-hand side by $t_j$ and
$t_j\frac{\partial}{\partial t_j}-s_j-1$, and on the right-hand side by $\tau_j$ and $-s_j$.
\end{lemme}

\begin{preuve}
It is sufficient to apply $p$ times the computation of the case B of the lemma \ref{DAnul-pvar}. Hence, the Koszul complex is isomorphic to its $0^{\mathrm{th}}$ cohomology group, \emph{i.e.} $\pi_* \mathcal{O}_{\C}^{an}$, \emph{via}  the application $\dsum_{\alpha \in (\N \setminus \{0\})^p}b_\alpha t^\alpha \flto b_{1,\dots,1}$.\\
Moreover, an easy checking shows that the actions of $t_j$ and
$t_j\frac{\partial}{\partial t_j}-s_j-1$ become respectively the actions of $\tau_j$ and $-s_j$ \emph{via} this morphism.  
\end{preuve}

Hence, we get the isomorphism of $\mathrm{D}$-modules
$$ \mathcal{K}[p] \cong \pi_* \mathcal{O}_{\C^p}$$
where $t_j$ and $t_j\partial_{t_j}$ act on the left-hand side by $t_j$ and
$t_j\frac{\partial}{\partial t_j}-s_j-1$, and on the right-hand side by $\tau_j$ and $-s_j$.

\subsubsection{The end of the proof of theorem \ref{th-mellin}}

Thanks to proposition \ref{prop-mellin-sol} and the preceeding paragraphs, we get
$$\mathfrak{M}(\mathcal{S}ol^{<0}(\mathcal{M}))  \cong \textbf{R}\mathcal{H}om_{\mathrm{D}} (  \mathrm{M} ,
\mathcal{K} )[p] \cong \textbf{R}\mathcal{H}om_{\mathrm{D}} (  \mathrm{M} ,
\pi_* \mathcal{O}_{\C^p} )$$ 
where $t_j$ and $t_j\partial_{t_j}$ act on $\pi_* \mathcal{O}_{\C^p}$ by $\tau_j$ and $-s_j$.\\

Hence, we get 
$$\mathfrak{M}(\mathcal{S}ol^{<0}(\mathcal{M}))  \cong
\textbf{R}\mathcal{H}om_{\C[s]\langle \tau,\tau^{-1}\rangle}
(\mathfrak{M}(\mathcal{M}) , \pi_* \mathcal{O}_{\C^p}) $$ 
dans $D^b(\mathfrak{Mod}(\mathcal{O}_{\C^p / \Z^p}))$.

\end{document}